\documentclass[a4paper]{amsart}\usepackage{amssymb}

\usepackage{verbatim}
\usepackage{amscd}
\def\vrt{{{\operatorname{vrt}}}}
\def\cnc{{{\operatorname{cnc}}}}
\def\trg{{{\operatorname{trg}}}}
\def\src{{{\operatorname{src}}}}
\def\per{{{\operatorname{per}}}}
\def\card{{{\operatorname{card}}}}

\def\id{{{\operatorname{id}}}}

\numberwithin{equation}{section}
\theoremstyle{plain}
\newtheorem{lemma}{Lemma}[section]
\newtheorem{theorem}{Theorem}[section]
\newtheorem{corollary}{Corollary}[section]
\newtheorem{proposition}{Proposition}[section]
\newtheorem*{Lemma A}{Lemma A}
\newtheorem*{Lemma B}{Lemma B}
\newtheorem*{Lemma C}{Lemma C}
\theoremstyle{definition}

\theoremstyle{remark}

\begin{document}
\title [On certain labelled directed graphs]{On certain labelled directed graphs \\ of symbolic dynamics}
\author{Wolfgang Krieger}
\begin{abstract}
We consider families of coded systems that contain the Dyck shifts and that are closed under topological conjugacy.
We introduce a notion of hyposynchronization of subshifts. 
We introduce a notion of restricted complexity of hyposynchronizing subshifts. 
Restricted complexity of hyposynchronizing subshifts is not invariant under topological conjugacy.
  We construct explicitly a family of of hyposynchronizing subshifts of restricted complexity that extends the family of Dyck shifts. The subshifts in this family are characterized by their restricted complexity in conjunction with a certain set of  invariants of topological conjugacy that includes the Artin-Mazur zeta function.
\end{abstract}
\maketitle

\section{Introduction}
Let the symbol $\Sigma$ denote a finite alphabet. By a subshift is meant a dynamical system that is obtained by restricting the left shift $S_\Sigma$ on $\Sigma^{\Bbb Z}$ to one of its closed invariant subsets. 
For an introduction to the theory 
of subshifts see \cite {Ki} or  \cite{LM}. A word in the symbols of $\Sigma$ is said to be admissible for a subshift $X \subset \Sigma^{\Bbb Z}$ if it appears somewhere in a point of $X$. We denote the language of admissible words of a subshift $X$ by $\mathcal L(X)$.

As a subshift $X \subset \Sigma^{\Bbb Z}$  is uniquely determined by $\mathcal L(X)$, it 
can be defined by specifying $\mathcal L(X)$. For this purpose  $\Sigma$-labelled directed graphs can be used: A $\Sigma$-labelled directed graph,  in which every vertex has an incoming edge and an outgoing edge defines the subshift that has as its language of admissible words the set of label sequences of the finite directed paths in the graph. We say that the graph presents the subshift \cite {Kr7}.

A $\Sigma$-labelled directed graph is said to be 1-right resolving if every vertex of the graph has for every symbol $\sigma \in \Sigma$ at most one outgoing edge that carries the label $\sigma$. A 1-right resolving $\Sigma$-labelled directed graph that presents the subshift $X \subset \Sigma^{\Bbb Z}$ can be viewed as a real time deterministic automaton with all vertices as initial and final states that 
recognizes the language $\mathcal L(X)$.

Arguably the most important equivalence relation for subshifts is topological conjugacy. This is the only equivalence relation for subshifts that we consider in this paper. Basic classes of subshifts that are closed under topological conjugacy are the subshifts of finite type 
\cite {P, W},  the sofic systems \cite {W} and the coded systems  \cite {BH3}. The subshifts of finite type belong to the family of sofic systems and irreducible sofic systems are coded. Families of coded systems that are closed under topological conjugacy can be described by means of known invariants of flow equivalence \cite{Kr4,HK1, HK2, Kr8, CS, M2}.

Prominent examples of coded systems that are not sofic are the Dyck shifts \cite {Kr1}. These serve as prototypical examples for the type of subshifts that we study in this paper. 
To recall the construction of the Dyck shifts,
let $N> 1$, and let  
$
\alpha^-(n), \alpha^+(n), 0 \leq n < N,
$
 be a set of generators of the Dyck inverse monoid (the polycyclic monoid \cite {NP}) $\mathcal D_N$, which is the graph inverse monoid of the directed graph with one vertex and N loops. The generators  satisfy the relations
$$
\alpha^-(n) \alpha^+(m) =
\begin{cases}
\bold{1}, &\text{if  $n = m$}, \\
0, &\text {if $n \neq m$}. \quad \quad 1\leq n,m \leq N.
\end{cases}
$$
The Dyck shift  $D_N$ is defined as the subshift
$$
D_N \subset( \{ \alpha^-(n): 0 \leq n < N \} \cup\{ \alpha^+(n):0 \leq n < N \})^\Bbb Z
$$
with admissible words $(\sigma_i)_{1 \leq i \leq I  } , I \in \Bbb N,$ of $D_N, N > 1,$ given by the condition 
\begin{align*}
\prod_{1 \leq i \leq I } \sigma_i \neq 0. 
\end{align*}

Neutral periodic points of the Dyck shifts appeared in \cite{HI} in connection with the problem of embedding subshifts of finite type into the Dyck shifts (and more general target shifts, see \cite{HIK}). 
A periodic point $p$ of period $\pi$ 
of the Dyck shift $D_N$ is neutral if there exists an index 
$i \in \Bbb Z$ such that
$$
\prod_{i \leq j < \pi}  p_j = \bold{1}.
$$
The definition of a neutral periodic point of a subshift was given in  \cite{Kr9}. 
In this paper we exploit the structure of a subshift that emanates from its set of neutral periodic points.

Without aiming for closure under topological conjugacy families of coded systems containing the Dyck shifts have been described:
By means of  $\mathcal R$-graph semigroups the family of $\mathcal R$-graph shifts was constructed in \cite{Kr9} that contains the Dyck shifts (see also \cite{HIK} for a family that includes the $\mathcal R$-graph shifts). In \cite{BH2} there was described a family of generalized Dyck shifts. In \cite{BBD2, BBD3} there was introduced a family of subshifts under the name of sofic Dyck shifts, together with the Dyck automata, that recognizes these,  and with a subfamily of finite type Dyck shifts \cite{BBD1}.
The subshifts that arise from 
$\mathcal R$-graph semigroups are finite type Dyck shifts \cite{Kr1}. The subshifts that arise from graph inverse semigroups are the Markov Dyck shifts \cite{KM1}. Following \cite{KM1} the case of Markov-Dyck shifts of rotationally homogeneous finite trees was considered in \cite{Kr10}.

In  \cite {Kr5} the notion of semisynchronization of subshifts was introduced. Semisynchronization of subshifts is an invariant of topological conjugacy. Semisynchronizing subshifts are coded and the Dyck shifts are semisynchronizing. In this paper we introduce a notion of hyposynchronization of subshifts as a companion notion for the notion of semisynchronization of subshifts. Hyposynchonization is an invariant of topological conjugacy. Hyposynchronizing subshifts are coded and Dyck shifts are hyposynchronizing. Other notions of symbolic dynamics that express a quality of synchronization, weaker than synchronization, and 
that are invariant under topological conjugacy, are $\lambda$-synchronization \cite{KM1} and Property D \cite{Kr6}.

In addition to the family of hyposynchronizing subshifts we consider in this paper a series of families of coded systems that are closed under topological conjugacy and that contain the Dyck shifts. 
We obtain these families by imposing invariant conditions on the coded systems. We organize these conditions in a list\footnote{We accept into this list only invariant conditions that can be formulated in such a way that it is possible, in principle, to verify them given any presentation of the subshift.}
that comprises Conditions $\bold H$01 - $\bold H $20. The conditions on this list are not meant to be independent.

We introduce the 
notion of a neutral point of a subshift
$X \subset \Sigma^\Bbb Z$: A point $x \in X$ is neutral if for every time point $i \in \Bbb Z$ there exists a time span $\ell (i)\in \Bbb N$ such that the future 
$(x_j)_{i < j < \infty}$
that is given by the time point is compatible with every past that in $X$ can precede the event $(x_j)_{i - \ell(i) < j < i}$,
and also such that the past $(x_j)_{- \infty <  j < i}$ that is given by the time point is compatible with every future that in $X$ can follow the event $(x_j)_{i  < j < i+\ell(i)}$.

We say that a point $x$ of a subshift 
$X \subset \Sigma^{\Bbb Z}$ is left (right) asymptotic to a neutral periodic point $p$ of $X$ if there exists an index $i_{\circ}\in \Bbb Z$ such that $x_i = p_{i_{\circ}+ i}$ for $i < 0$ ($i>0$).
We say that an orbit $\frak x$ of a subshift 
$X \subset \Sigma^{\Bbb Z}$ is left (right) asymptotic to a neutral periodic orbit $\frak p$ of $X$ if there exists $x\in \frak x$ and $p\in \frak p$,
such that $x$ is left (right) asymptotic to $p$.

The contents of this paper is as follows:
In a preliminary Section 2 we introduce notation and terminology for subshifts and for  directed graphs. We also recall prerequisites. 

In Section 3 we list conditions $\bold H$01 - $\bold H$11 on a subshift $X\subset \Sigma^\Bbb Z$. Conditions
$\bold H01$ and $\bold H$02 are irreducibility conditions that concern the orbits of a subshift $X\in \Sigma^\Bbb Z$ that are left asymptotic to a neutral periodic orbit of $X$ and also right asymptotic to a neutral periodic orbit of $X$. Condition $\bold H$03  allows to construct an operation that connects suitably structured points and orbits  of $X$. We call this operation a concatenation. The purpose of
$\bold H$04 - $\bold H$10 is to associate invariantly to  $X$ an $\mathcal R$-graph
$G_{\mathcal R(X)}(\frak V(X), \mathcal E^-(X) \cup \mathcal E^+(X) )$(\cite {Kr9}, see Section 2).
Condition $\bold H$04 concerns the introduction into the set of neutral periodic points of  $X$ of a basic equivalence relation $\sim$: Neutral periodic points 
$q$ and $r$ of $X$ are $\sim$-equivalent if there exist  
 neutral points $u$ and $v$ of $X$ such that  $u$ is left asymptotic to $q$ and right asymptotic to $r$, and such that $v$ is left asymptotic to $r$ and right asymptotic to $q$. The vertex set $\frak V(X)$ is the set of $\sim$-equivalence classes. 
 Condition $\bold H$11 is a condition on the $\mathcal R$-graph
$G_{\mathcal R(X)}(\frak V(X), \mathcal E^-(X) \cup \mathcal E^+(X) )$ that implies that the subshift $X$ is not synchronizing.
 The  
$\mathcal R$-graph shifts are prototypical examples of subshifts that satisfy Conditions  
$\bold H$01 - $\bold H$11. 

In Section 4 we list Conditions $\bold H$12 - $\bold H$18. Condition $\bold H$12 yields the defining property for the introduction of the notion of a hyposynchronizing subshift.
Hyposynchronizing subshifts are the subshifts that satisfy Conditions $\bold H$01 - $\bold H$12.
By means of Condtions $\bold H$12 - $\bold H$16 we obtain for 
$\textcircled{\tiny V}\in \frak V(X)$ a 1-right resolving strongly connected labelled directed graph $G_{\textcircled{\tiny V}}^{hypo}(X)$ that presents the hyposynchronizing subshift $X$.
Conditions  $\bold H$17 and $\bold H$18 concern the cycles in the graphs 
$G_{{\textcircled{\tiny V}}}^{hypo}(X)$. The two conditions allow to
assign a height to the vertices of $G_{\textcircled{\tiny V}}^{hypo}(X)$. 
For the Dyck shifts this height coincides with the natural one.

Taking the $\mathcal R$-graph shifts as motivating examples we say that a hyposynchronizing subshift has restricted complexity if it can be presented by a 1-right resolving pushdown automaton in which the stack content of a configuration determines the internal state of the configuration\footnote{This means that internal states do not come up in the description of the pushdown automata. Such automata are called "simple" in \cite[p.139]{ABB}. We use the term "complexity". The idea is that the number of internal states of a pushdown automaton is a measure of its complexity.}.
In Section 5 we describe the explicit construction of a family $\mathcal Y$ of hyposynchronizing subshifts of restricted complexity. The family $\mathcal Y$ contains the Dyck shifts and all subshifts in $\mathcal Y$ are patterned after the Dyck shifts. The familiy 
$\mathcal Y$ intersects the family of generalized Dyck shifts that was introduced in  \cite {BH2}. In this section we list $\bold H$19 and $\bold H$20.   
The subshifts in the family 
$\mathcal Y$ are characterized by satisfying $\bold H$19 and $\bold H$20 and  by their Artin-Mazur zeta function. 

Section 6 contains examples.
\section{Preliminaries}

We introduce notation and terminology for subshifts  $X \subset \Sigma^\Bbb Z$.
We denote the set of periodic points of $X$ by $P(X)$.
The period of a point $p\in P(X)$ we denote by $\per(p)$.
 The set of orbits of $X$ we denote by $\mathcal O(X)$.
 
Given a subshift $X \subset \Sigma^\Bbb Z$, and $I^-, I^+ \in \Bbb Z, I^-\leq I^+$
we call
$(\sigma_i)_{I^-\leq i \leq I^+}\in \Sigma^{[I^-, I^+]}$ a block
and we use the notation
$$
x_{[I^-, I^+]} = (x_i)_{I^- \leq i \leq I^+},\quad \quad \quad    x \in X.
$$
For  a subset $A$ of $X$ we use the notation
$$
A_{[I^-, I^+]}   = \{x_{[I^-, I^+]}  : x \in A\}.
$$
We use similar notations in the case that indices range in half-open, left infinite or right infinite intervals.

We identify a block in $\Sigma^{[1, n]}, n \in\Bbb N$ with the word that it carries.
More generally use the same notation for a block in $ \Sigma^{[I^-, I^+]}$ and for the word that the block carries. From the context it is always clear what is meant, the block or the word.

We use notations like $x^{(-\infty, i]}$($x$) for the left (right) infinite points in $X_{(-\infty,i]}$($ X_{[i,\infty)}$).

We set
$$
\Gamma^{\langle + \rangle}(a) = \{ x^{(j, \infty)} \in X_{(j, \infty)}: a x^{(j, \infty)}
 \in X_{[i, \infty)}\}, \quad 
a \in  x_{[i,j]}, \quad i,j \in \Bbb Z, i \leq j.
$$
The notation $\Gamma^{\langle - \rangle}$ has the symmetric meaning.
We also set
$$
\omega^{\langle + \rangle}(a) = \bigcap_{x^{(\-\infty,}\in \Gamma^{\langle - \rangle}(a)} \{ x^{(j, \infty)} \in \Gamma^{\langle + \rangle}(a):x^{(-\infty, i)}a x^{(j, \infty)} \in X\}, 
$$
$$
 a \in  x_{[i,j]}, \quad i,j \in \Bbb Z, i \leq j.
$$
The notation $\omega^{\langle - \rangle}$ has the symmetric meaning.

Consider subshifts $X \subset \Sigma^\Bbb Z, \bar X \subset \bar \Sigma^\Bbb Z $
and a continuous shift commuting surjection
$$
\varphi: X \to \bar X.
$$
The mapping $\varphi$ is implemented by a block map \cite{Ki, LM}. For such an implementation there are given $I^-, I^+ \in \Bbb Z, I^-\leq I^+$. We refer to $I^-$ as a memory and to $I^+$ as an anticipation. The block map  $\Phi$ is a map
$$
\Phi: \mathcal L_{ I^+ - I^-  +1} \to \bar \Sigma.
$$
The implemetation of $\varphi$ by $\Phi$ is given by
$$
\varphi((x_i)_{ -\infty < i < \infty} ) = 
(\Phi ((x_i)_{I^- + j \leq i \leq I^+ + j} ))_{ -\infty < j < \infty}, \quad \quad x\in X.
$$
In the case that $\varphi$ is a conjugacy and the memory equals the anticipation the implementation renames the symbols of the alphabet $\Sigma$. We call such a conjugacy a literal conjugacy.

The following lemma can be used in invariance proofs.

\begin{lemma}
Let there be given subshifts 
$X \subset \Sigma^\Bbb Z, \widetilde{ X} \subset \widetilde{ \Sigma}^\Bbb Z,   $
and a topological conjugacy
$$
\varphi :   \widetilde{ X} \to X.
$$
Let $\varphi$ be implemented by a block map $\Phi$ and let $\varphi^{-1}$ be implemented by a block map $\widetilde{\Phi}$. Let $[-L, L ]$ be a coding window for both, 
 $\Phi$ and $\widetilde{\Phi}$.
 
 Let $\widetilde{x } \in \widetilde{X}  $, and set
 \begin{align}
 x = \varphi(\widetilde{x }) .
  \end{align}
  Let
   \begin{align}
   K > 2L.
    \end{align}
    and let
    \begin{align}
    x_{[1, \infty)} \in \omega_\infty^{\langle + \rangle}(x_{(-K, 0]}).
     \end{align}
    Then
    \begin{align}
    \widetilde{x }_{(L, \infty)} \in  \omega_\infty^{\langle + \rangle}(\widetilde{x}_{(-K-L, L]}).
     \end{align}
\end{lemma}
\begin{proof}
The task is to show that 
 \begin{align}
 ( \widetilde{y}_i)_{-\infty < i < K - L} \in 
 \Gamma_\infty^{\langle - \rangle}(\widetilde{x}_i)_{-K- L \leq i \leq L},
  \end{align}
implies that
\begin{align}
\widetilde{z} =((  \widetilde{y}_i)_{-\infty < i < -K-L } ,( \widetilde{x}_i)_{-K- L\leq i < \infty}) \in \widetilde{X}.
\end{align}
We construct a point $z\in X$ such that
\begin{align}
\widetilde{z} =\varphi^{-1}(z).
\end{align}
By (2.7)
\begin{align}
((\widetilde{y}_i)_{-\infty < i < -K-L},  (\widetilde{x}_i)_{-K- L\leq i \leq L}) \in   \widetilde{X}_{(-\infty, L])}.
\end{align}
Therefore
\begin{align}
\Phi( (\widetilde{y}_i)_{-\infty < i <K- L} , 
(\widetilde{x}_i)_{- K- L\leq i \leq L}) \in
 {X}_{(-\infty, 0]}.
\end{align}
By (2.3) and (2.4)
\begin{align}
\Phi((\widetilde{x}_i)_{-K- L\leq i < L} )  = (x_i)_{-K\leq i \leq 0}.
\end{align}
Set
\begin{align}
(y_i)_{-\infty < i < -K} = \Phi(( \widetilde{y}_i)_{-\infty < i < -K+L}).
\end{align}
By (2.11) and  (2.12)
$$
(y_i)_{-\infty < i \leq -K}\in \Gamma_\infty^{\langle - \rangle}(({x}_i)_{-K < i \leq 0}).
$$
By (2.5)
$$
z = ((y_i)_{-\infty < i \leq -K}, (y_i)_{-K \leq i < \infty})\in X.
$$

By (2.3) and (2.4)
$$
\varphi^{-1} (z)= (\widetilde{\Phi}(z_{[i-L, i + L]})_{-\infty < i < \infty},
$$
and (2.9)  follows from  (2.3), (2.4) and  (2.13)       .     
\end{proof}

We introduce notation and terminology for labelled directed graphs. We denote a finite directed graph  with
vertex set ${\frak V}$ and edge set 
${\mathcal E}$ by $G(\frak V, \mathcal E)$.
As notation for the source vertex and target vertex of an edge or path in a directed graph we use $\src$ and $\trg$.

A $\Sigma$-labelled directed graph is called 1-right resolving  if  for every vertex $V$  and every label $\sigma\in \Sigma$ there is at most one edge that leaves $V$ and that carrries the label $\sigma$. We denote the target vertex of the edge that leaves the vertex $V$ and that carries the label $\sigma$ by $\tau_\sigma(V)$.
A bi-infinite directed path in the 1-right resolving $\Sigma$-labelled directed graph
$G(\frak V, \mathcal E)$ can be written as a bi-infinite sequence
$$
(V_i, x_I)_{i\in \Bbb Z} \in ( \frak V \times\Sigma)^\Bbb Z,
$$
such that 
$$
\tau_{x_i}(V_i) = V_{i+1},              \quad \quad \quad i \in \Bbb Z.
$$
The edge shift $E(G(\frak V, \mathcal E))$ of the 1-right resolving $\Sigma$-labelled directed graph 
$G(\frak V, \mathcal E)$ is the set of bi-infinite directed paths in the graph that is equipped with the left shift. This edge shift can be viewed as a countable state topological Markov shift.

We recall from \cite{Kr9} the notion of  
an $\mathcal R$-graph.
Let there be given 
a finite  directed graph $G(\frak V, \mathcal E)$. Assume also given a partition 
$$
\mathcal E = \mathcal E^-  \cup\mathcal E^+.
$$
We set
\begin{align*}
& \mathcal E^- (\frak q,\frak r) = \{ e^- \in  \mathcal E^- : s(e^-) = \frak q,\  t( e^-) =  \frak r \},
\\
& \mathcal E^+(\frak q,\frak r) = \{ e^- \in  \mathcal E^+ : s(e^+) = \frak r,\  t( e^+) =  \frak q \}, \qquad  \frak q,\frak r \in \frak V.
\end{align*}
 We assume that $ \mathcal E^- (\frak q,\frak r)  \neq \emptyset$ if and only if $  \mathcal E^+(\frak q,\frak r) \neq \emptyset,  \frak q,\frak r \in \frak V$, and we assume that 
the  directed graph $G(\frak V  ,  \mathcal E^-   )$ is strongly connected, or, equivalently, that 
the  directed graph $G(\frak V  ,  \mathcal E^+ )$ is strongly connected.
Let there further be given complete heterogeneous relations 
$$
\mathcal R    (\frak q,\frak r) \subset   \mathcal E^- (\frak q,\frak r)   \times   \mathcal E^+(\frak q,\frak r) , \qquad \frak q,\frak r \in \frak V,
$$
and set
$$
\mathcal R = \bigcup_{ \frak q,\frak r \in \frak P} \mathcal R    (\frak q,\frak r) .
$$
The resulting structure, that we call an $\mathcal R$-graph, we denote by 
$G_\mathcal R(\frak P, \mathcal E^-\cup\mathcal E^+)$. 

We also recall 
the construction from an $\mathcal R$-graph $G_\mathcal R(\frak V, \mathcal E^-\cup\mathcal E^+)$ of a semigroup (with zero)  $\mathcal S(G_\mathcal R(\frak V,   \mathcal E^-\cup \mathcal E^+   ))$ 
 as described in \cite {Kr9}.   
The semigroup $\mathcal S(G_\mathcal R(\frak V,   \mathcal E^-\cup \mathcal E^+))$ contains idempotents $\bold 1_{\frak p}, \frak p \in \frak V,$ and  has $\mathcal E$ as a generating set.
Besides $\bold 1_{\frak p}^2 = \bold 1_{\frak p},\frak p \in \frak P$, the defining relations are:
$$
f^- g^+ =\bold 1_{\frak q}, \quad f^- \in  \mathcal E^- (\frak q,\frak r), g^+ \in  \mathcal E^+(\frak q,\frak r) , (  f^- ,  g^+ ) \in \mathcal R    (\frak q,\frak r), \quad \frak q,\frak r \in\frak P,
$$
and
\begin{align*}
&\bold 1_{\frak q} e^- = e^- \bold 1_{\frak r} = e^-, \quad e^- \in  \mathcal E^- (\frak q,\frak r), \\
&\bold 1_{\frak r} e^+ = e^+ \bold 1_{\frak q} = e^+, \quad e^+ \in  \mathcal E^+ (\frak q,\frak r), \quad  \frak q,\frak r \in\frak V,
\end{align*}
$$
f^- g^+  = \begin{cases}{ \bold 1}_{{\frak q}}, &\text {if $ (  f^- ,  g^+ ) \in{ \mathcal R}({\frak q},{\frak r}) $,}\\
0, & \text{if $ (  f^- ,  g^+ )\notin {\mathcal R}({\frak q},{\frak r}), \quad f^- \in  {\mathcal E}^- ({\frak q},{\frak r}), g^+ \in { \mathcal E}^+({\frak q},{\frak r)},        
\
{\frak q},  {\frak r} \in{\frak V},
 $}
\end{cases}
$$
and
$$
\bold 1_{\frak q}\bold 1_{\frak r}= 0,    \quad  \frak q,\frak r \in \frak P,\frak q \neq\frak r .
$$
We call $\mathcal S_\mathcal R(G(\frak V,   \mathcal E^-\cup \mathcal E^+ ))$ an $ \mathcal R$-graph semigroup. 
We write $\mathcal S^-(\frak  V,   \mathcal E^- )$($\mathcal S^+(\frak P,   \mathcal E^+ )$) for the set of non-zero elements of the subsemigroup of 
$\mathcal S(G_\mathcal R(\frak  V,   \mathcal E^-\cup \mathcal E^+ ))$, that is generated by 
$ \mathcal E^-$
($\mathcal E^+$).
For the case that the $\mathcal R$-graph has a single vertex compare \cite{Ke2, LDT, MN} \cite[Section 1.5] {R}\cite[Example 1.2]{Li}.

Special cases are the graph inverse semigroups of  finite directed graphs 
 $ G (\frak V,\mathcal E)$ (\cite {AH},\cite [Section 10.7.]{ La}). With the edge set 
 $\mathcal E^- = \{  e^-: e_\circ \in \mathcal E_\circ\}$ of a copy of 
 $G(\frak V,\mathcal E)$, and with the  edge set $\mathcal E^+ =
  \{  e^-: e \in \mathcal E\}$ of the reversal of $G(\frak V,\mathcal E)$, the 
 graph inverse semigroup 
 $\mathcal S( G(\frak V,\mathcal E))$ of $G(\frak V,\mathcal E)$
  is  the $\mathcal R$-graph semigroup of the partitioned graph 
 $G(\frak V, \mathcal E^-  \cup\mathcal  E^+ )  $ with the relations
$$
\mathcal R(\frak q, \frak r) = \{(e^- , e^+ ) : e  \in \mathcal E   , s(e ) = \frak q, t(e ) = \frak r \}, \quad   \frak q,\frak r \in \frak V.
$$

The  $\mathcal R$-graph shift 
$X(G_\mathcal R(\frak V, \mathcal E^- \cup \mathcal E^+))$ of the $\mathcal R$-graph 
$G_\mathcal R(\frak V, \mathcal E^- \cup \mathcal E^+)$ is the subshift with alphabet
$\mathcal E^- \cup \mathcal E^+$,
and admissible words 
 given by the condition 
\begin{align*}
\prod_{1 \leq i \leq I } \sigma_i \neq 0. 
\end{align*}
To avoid the case that the shift $X(G_\mathcal R(\frak V, \mathcal E^- \cup \mathcal E^+))$ is finite, assume that the graph $G(\frak V, \mathcal E^-)$ is not a circle.

 \section{Neutral periodic orbits}
  
 \subsection{Neutral periodic points and orbits}
 
 We consider a subshift $X \subset \Sigma^{\Bbb Z}$. We set
 \begin{align*}
 &J_i^{\langle -\rangle}(x) = \{ j \in \Bbb N: x_{[i,\infty)}\in 
 \omega ^{\langle +\rangle}(x_{[i-j,i)} )\},
 \\
 &J_i^{\langle +\rangle}(x) = \{ j \in \Bbb N: x_{[i,\infty)}\in 
 \omega ^{\langle +\rangle}(x_{[i-j,i)} )\}, \quad i \in \Bbb Z, x \in X.
\end{align*}
 We set
 $$
 A^{(-)}(x) = \bigcap _{i\in \Bbb Z} \{x\in X:J_i^{\langle -\rangle}(x) \not =\emptyset \},
 $$
 $$
 A^{(+)}(x) = \bigcap _{i\in \Bbb Z} \{x\in X:J_i^{\langle +\rangle}(x) \not =\emptyset \},
 $$
 $$
 A^0(x)= A^{(-)}(x) \cap A^{(+)}(x), 
 $$
 $$
 A^-(x)=A^{(-)}(x) \setminus A^0(x), \quad  A^+(x)=A^{(+)}(x) \setminus A^0(x).
 $$
 By Lemma 2.1 the shift invariant sets $A^-(x), A^0(x),A^+(x)$ are invariantly associated to the subshift $X$. We set
 $$
 P^0(X)= P(X) \cap A^0(x),
 $$
 $$
 \mathcal O^{(0)}(X) = \{\frak p \in \mathcal O(X): \frak p \subset  A^0(X)\}.
 $$
The points in $P^0(X)$ are the neutral periodic points of $X$. The orbits in $\mathcal O^0(X)$ are the neutral periodic orbits of $X$.

 \smallskip
 
We denote by
$\mathcal F(\frak p , \bullet) $($\mathcal F(\bullet,\frak p ) $)
the set of orbits in $\mathcal O(X) \setminus \{\frak p\}$ that are  left (right) asymptotic to
$\frak p\in \mathcal P^0(X)$.
 For $\frak p \in \mathcal P^0(X) $ and orbits 
 $$
 \frak x^{\langle - \rangle}\in \mathcal F(\bullet, \frak p)
  , \quad \frak x^{\langle + \rangle} \in \mathcal F(\frak p, \bullet),
 $$
 and points
  $x^{\langle - \rangle} \in  \frak x^{\langle - \rangle} $, $x^{\langle + \rangle} \in  \frak x^{\langle + \rangle},$ we set
$$ 
 I^+(x^{\langle - \rangle}) =\max ( \{  i\in \Bbb Z: x_i = x_{i - \card(\frak p)} \}, \quad
  I^-(x^{\langle + \rangle}) =\min  \{  i\in \Bbb Z: x_i = x_{i + \card(\frak p)} \}.   
 $$

 We set
 $$
 \mathcal F(\frak q , \frak r) = \mathcal F(\frak  q,\bullet ) \cap\mathcal F(\bullet,\frak  r ),
 \quad  \quad\frak q,\frak r \in \mathcal O^0(X),
 $$
$$
  \mathcal F(X) = \bigcup_{\frak q ,\frak r \in   \mathcal O^{(0)}(X) }
   \mathcal F(\frak q , \frak r).  
  $$
  We also set
  \begin{align*}
 &\mathcal F^{0}(\frak  p,\bullet )
 = \{ \frak x \in \mathcal F(\frak  p,\bullet )): 
 \frak x \subset A^0(X) \},  
 \\  
 & \mathcal F^{0}(\bullet,\frak  p ) =
  \{ \frak x \in \mathcal F(\bullet,\frak  p )): 
  \frak x \subset A^0(X) \},
   \quad \quad \quad \frak p \in \mathcal O^0(X) .
    \end{align*}
    
 As irreducibility conditions to be imposed  on the subshift $X \subset \Sigma^{\Bbb Z}$ we adopt the following two conditions $\bold H01$ and $\bold H02$:

 \smallskip

\smallskip
 \noindent
 $\bold H01$ For $\frak q,\frak r \in \mathcal O^0(X)$ the sets 
 $\mathcal F^{-}(\frak  q,\frak r ) $ and $\mathcal F^{+}(\frak  q,\frak r )$ are not empty.

 \smallskip
 \noindent
 $\bold H02$ For $\frak p \in \mathcal O^0(X)$ the set 
 $\mathcal F^{0}(\frak p,\frak p)$ is not empty.
 
 \smallskip
 
  \subsection{Concatenation of orbits}
 
Let there be given an orbit $\frak p \in \mathcal P^0(X) $ and orbits 
$$
 \frak x^{\langle - \rangle}\in \mathcal F(\bullet, \frak p)
  , \quad \frak x^{\langle + \rangle} \in \mathcal F(\frak p, \bullet).
$$
Let points
$$
x^{\langle - \rangle} \in  \frak x^{\langle - \rangle}, \quad x^{\langle + \rangle} \in  \frak x^{\langle + \rangle},
 $$
 be given by the condition
 $$
 I^+ (x^{\langle - \rangle}) = I^-(x^{\langle + \rangle}) = 0,
 $$
 and let points $p^{\langle - \rangle}, p^{\langle + \rangle} \in \frak p$ be given by the condition that $x^{\langle - \rangle}$ is right asymptotic to $p^{\langle - \rangle}$ and
 $x^{\langle + \rangle}$ is left asymptotic to $p^{\langle + \rangle}$.
 Let $D(p^{\langle - \rangle}, p^{\langle + \rangle} )\in  [0, \thinspace \card( \frak p)) $ be given by the condition
 $$
 S_X^{D(p^{\langle - \rangle}, p^{\langle + \rangle})}p^{\langle - \rangle} = 
 p^{\langle + \rangle}.
 $$
 For $m \in \Bbb N$ we construct points $z^{\langle -, m \rangle}$ and 
 $z^{\langle +, m \rangle}$ in $\Sigma^\Bbb Z$ that have the same orbit under 
 $S_\Sigma $, and that are symmetric to each another:
  
 \begin{align*}
 z^{\langle -, m \rangle}_{(-\infty, -m \times \card(\frak p) - 
 D(p^{\langle - \rangle}, p^{\langle + \rangle}))} &=
 x^{\langle + \rangle}_{(-\infty, I^{\langle + \rangle}(x^{\langle - \rangle})) },
 \\
  z^{\langle -, m \rangle}
  _{[- m \times \card({\frak p}) -  D(p^{\langle - \rangle}, p^{\langle + \rangle})), 0]}&=
  p^{\langle - \rangle} _{[- m \times \card({\frak p}) -  D(p^{\langle - \rangle}, p^{\langle + \rangle})), 0]},
 \\
 z^{\langle -, m \rangle}_{(0, \infty )}&=
 x^{\langle + \rangle}_{(I^{\langle - \rangle}(x^{\langle + \rangle}, \infty)},
\end{align*}

  \smallskip
  
   \begin{align*}
  z^{\langle +, m \rangle}_{(-\infty, 0 )}&=
 x^{\langle - \rangle}_{(-\infty, I^{\langle + \rangle}(x^{\langle - \rangle}))}, 
 \\
 z^{\langle +, m \rangle}_{[0, m \times \card(\frak p) + D(p^{\langle - \rangle}, p^{\langle + \rangle})]} &=
 p^{\langle + \rangle}_{[0, m \times \card(\frak p) + D(p^{\langle - \rangle}, p^{\langle + \rangle})]},
 \\
  z^{\langle +, m \rangle}
  _{(m \times \card({\frak p}) +  D(p^{\langle - \rangle}, p^{\langle + \rangle})), \infty)}&=
  x^{\langle + \rangle}_{( I^{\langle - \rangle}(x^{\langle + \rangle}), \infty) },
 \end{align*}

 \bigskip
  
We denote the orbit of $(\Sigma^\Bbb Z, S_\Sigma)$ that contains the points 
$z^{\langle -, m \rangle}$  and $z^{\langle +, m \rangle}$ by
$$
 \cnc_m \langle \frak x^{\langle - \rangle},\frak x^{\langle + \rangle} \rangle,
$$
and we call this orbit a concatenation of the orbits $\frak x^{\langle -\rangle}$ and 
$\frak x^{\langle + \rangle}$. For $x^{\langle - \rangle} \in   \frak x^{\langle - \rangle},
x^{\langle + \rangle} \in   \frak x^{\langle + \rangle},$ we say that the points in
$ \cnc_m \langle \frak x^{\langle - \rangle},\frak x^{\langle + \rangle} \rangle$ are concatenations  of  $x^{\langle - \rangle}$ and  
$x^{\langle + \rangle}.$ Note that the orbits $\cnc_m \langle \frak x^{\langle - \rangle},\frak x^{\langle +, \rangle} \rangle, m\in \Bbb N$, are distinct.
  
 As a condition to be imposed on a subshift $X \subset \Sigma^{\Bbb Z}$ that satisfies 
 $\bold H01$ and $\bold H02$ we adopt the condition
  
 \smallskip
 \noindent
 $\bold H03$.
 There exists $M \in \Bbb N$ such that the following holds for $m> M$:
 \noindent
 For $\frak p\in \mathcal P^0(X)$ and for
 $$
 \frak x^{\langle - \rangle} \in \mathcal F( \bullet, \frak p), \quad \quad 
 \frak x^{\langle + \rangle} \in \mathcal F(\frak p, \bullet),
 $$
 one has that the orbit 
 $\cnc_m\langle \frak x^{\langle - \rangle}  , \frak x^{\langle +\rangle} \rangle$
 belongs to $ \mathcal O(X)$
 if and only if the orbit
 $\cnc_M\langle \frak x^{\langle - \rangle}  , \frak x^{\langle +\rangle} \rangle$
 belongs to $ \mathcal O(X)$.
 
 \smallskip

 \begin{proposition}
 $\bold H03$ is an invariant of topological conjugacy.
 \end{proposition}
 \begin{proof}
 Let there be given subshifts $X \subset \Sigma^\Bbb Z$, $\bar{X} 
 \subset \bar\Sigma^\Bbb Z$, and a topological conjugacy
 $$
 \psi: \bar{X} \to X
 $$
 that is implemented by a block map
 $$
 \Psi: \bar{X}_{[-L, L]} \to \Sigma.
 $$
 We assume that $X$ satisfies  $\bold H03$ and we prove that $\bar{X}$ satisfies 
 $\bold H03$ with parameter
 $$
 \bar M = M +2L.
 $$
 
 Let $\bar{\frak p} \in \mathcal O(\bar{X})$, and let there be given orbits
 $$
 \bar{\frak x}^{\langle - \rangle} \in \mathcal F(\bullet , \bar{\frak p} ), \quad \quad
 \bar{\frak x}^{\langle +\rangle} \in \mathcal F( \bar{\frak p}, \bullet).
 $$
 Set
 $$
 {\frak x}^{\langle - \rangle} = \psi( \bar{\frak x}^{\langle - \rangle}), \quad \quad 
 {\frak x}^{\langle + \rangle} = \psi( \bar{\frak x}^{\langle + \rangle}).
 $$
 Let $\bar{m} \geq \bar M$. There exists $m \geq \bar{m}$, such that
 $$
 \psi(\cnc_{\bar{m}} 
 \langle \bar{\frak x}^{\langle - \rangle} , \bar{\frak x}^{\langle -+\rangle} \rangle ) =
 \cnc_{{m}} 
 \langle {\frak x}^{\langle - \rangle} , {\frak x}^{\langle +\rangle} \rangle. 
 $$
 Using the assumption that $X$ satisfies $\bold H03$ we have
 $$
 \cnc_{\bar{m}} 
 \langle \bar{\frak x}^{\langle - \rangle} , \bar{\frak x}^{\langle +\rangle} \rangle  
 \in \mathcal O(\bar{X}),
 $$
 if and only if
 $$
 \cnc_{{m}} 
 \langle {\frak x}^{\langle - \rangle} , {\frak x}^{\langle +\rangle} \rangle \in \mathcal O({X}),
 $$
 if and only if
 $$
  \cnc_{{M}} 
 \langle {\frak x}^{\langle - \rangle} , {\frak x}^{\langle +\rangle} \rangle \in \mathcal O({X}).
 $$
 It follows for $ \bar{m}>\bar{M}$ that
 $$
\cnc_{\bar{m}} 
 \langle \bar{\frak x}^{\langle - \rangle} , \bar{\frak x}^{\langle +\rangle} \rangle 
 \in \mathcal O(\bar{X}),
 $$
 if and only if
 $$
 \cnc_{{M}} 
 \langle {\frak x}^{\langle - \rangle} , {\frak x}^{\langle +\rangle} \rangle \in\mathcal O({X}),
 $$
 if and only if
 $$
\cnc_{\bar{M}} 
 \langle \bar{\frak x}^{\langle - \rangle} , \bar{\frak x}^{\langle +\rangle} \rangle  \in \mathcal O(\bar{X}). \qed
 $$
 \renewcommand{\qedsymbol}{}
  \end{proof}
 
 We introduce an equivalence relation $\approx$ into the set $\mathcal F(X)$:
Given orbits $\frak q. \widetilde{\frak q},\frak r, \widetilde{\frak r}\in \mathcal O^0(X),$
we say that  orbits 
$$
\frak x\frak \in \mathcal F(\frak q , \frak r), \qquad
\widetilde{\frak  x}\in \mathcal F( \widetilde{\frak q},\widetilde{\frak r} ),
$$
are 
$\approx$-equivalent, if
$$
\frak q = \widetilde{\frak q},  \qquad \frak r = \widetilde{\frak r}, 
$$ 
and if one has for all choices of
$$
 \frak u\in \mathcal F(\bullet , \frak q) , \qquad  \frak v\in \mathcal F(\frak r, \bullet) ,
$$
that
$$
\cnc_m \langle \cnc_m \langle \frak u, \frak x \rangle , \frak v \rangle \in \mathcal O(X),
$$
if and only if
$$
\cnc_m \langle \cnc_m \langle \frak u, \widetilde{\frak  x} \rangle , \frak v \rangle \in \mathcal O(X).
$$
One checks that $\approx$ is an equivalence relation.
  
More generally, 
consider a subshift $X \subset \Sigma^{\Bbb Z}$ that satisfies $\bold H01$ - $\bold H02$ and  $\bold H03$ with parameter $M$, and let there be given $K \in \Bbb N$, orbits 
 $$
 \frak p^{(k)} \in \mathcal O^0(X), \quad \quad 0\leq k \leq K,
 $$
 and orbits
 $$
 \frak x^{(k)}  \in \mathcal F(\frak p^{(k-1)}  , \frak p^{(k)}  ), \quad \quad 1\leq k \leq K.
 $$
 Let there also be given
 \begin{align}
 m_k \geq M, \quad \quad 1\leq k < K. 
 \end{align}
 We obtain orbits
 $$
 \frak y^{(k)},\quad \quad 1\leq k \leq K, 
 $$
 by 
  \begin{align*}
 &\frak y^{(1)} =  \frak x^{(1)},
  \\
  & \frak y^{(k)}=
   \cnc_{m_k}\langle \frak y^{(k-1)},\frak x^{(k)}  \rangle,
   \quad \quad 1< k \leq K. 
   \end{align*}
   For the orbit $\frak y^{(K)}$ we use the notation 
   \begin{align*}
    \cnc \langle (\frak x^{(k)})
   _{1\leq k \leq K}  \rangle.
   \end{align*}
   A more specific notation for this orbit would be
   \begin{align*}
   \cnc_{(m_k)_{1 \leq k \leq K}}\langle (\frak x^{(k)})_{1\leq k \leq K}  \rangle.
   \end{align*}
However,  were there given another set  $(m^\prime_k)_{1 \leq k \leq K} $ of parameters, then one would have by $\bold H03$ that
$$
\cnc_{(m_k)_{1 \leq k \leq K}}\langle (\frak x^{(k)})_{1\leq k \leq K}  \rangle \approx
\cnc_{(m^\prime_k)_{1 \leq k \leq K}}\langle (\frak x^{(k)})_{1\leq k \leq K}  \rangle.
$$
Consequently the orbits on both sides of this equivalence should be considered as being on the same footing. We therefore suppress the set of parameters (3.1) in the notation with the understanding that any set can be used, provided that within a specific context, e.g. within a proof, this set of parameters is not changed.

\smallskip

\subsection{The associated $\mathcal R$-graph}
 
We consider a subshift $X \subset \Sigma^{\Bbb Z}$, that satisfies  $\bold H01 - \bold H03$.
 
 As a condition to be imposed on the subshift $X \subset \Sigma^{\Bbb Z}$ that satisfies  $\bold H01 - \bold H03$. we adopt 
 
 \smallskip
  
 \noindent
 $\bold H04$. For $\frak q, \frak r \in \mathcal O^0(X), 
 \frak q   \neq  \frak r,$ the set 
 $\mathcal F^0(\frak q, \frak r)$ is not empty if and only if the set
 $ \mathcal F^0(\frak r, \frak q)$ is not empty.
 
 \noindent
 
 We say that $\frak q, \frak r \in \mathcal O^0(X)$ are $\sim$-equivalent if $\mathcal F^0(\frak q, \frak r) $ is not empty.
 
 \begin{lemma}
 The $\sim$-equivalence of neutral periodic orbits of $X$ is an equivalence relation.
 \end{lemma}
\begin{proof} 
By $\bold H04$ $\sim$-equivalence is reflexive, by $\bold H02$ it is symmetric. It is also transitive: let $\frak q, \frak p, \frak r \in \mathcal O^0(X)$, and let
$$
\frak x \in \ \mathcal F^0(\frak q, \frak p) , \qquad \frak y \in \mathcal F^0(\frak p, \frak r) .
$$
Then also
$$
\cnc \langle \frak x ,\frak y   \rangle  \in \mathcal F^0(\frak q, \frak r). \qed
$$
\renewcommand{\qedsymbol}{}
\end{proof}

We denote the set of $\sim$-equivalence classes of neutral periodic orbits of $X$ by 
$\frak V(X)$. 

As a condition to be imposed on the subshift $X \subset \Sigma^{\Bbb Z}$ that satisfies  
$\bold H01 - \bold H04$. we adopt the following two conditions.

\smallskip

\noindent
$\bold H05$. The set $\frak V(X)$ is finite.

\smallskip 

\noindent
$\bold H06$. For $V \in \frak V(X)$ the set $\bigcup_{\frak p \in V} \frak p$ is dense in $X$.

 \smallskip

We introduce an equivalence relation $\cong$ into the set $\mathcal F(X)$. We do this in two stages. In the first stage we introduce an equivalence relation $\approx$ into the set
$\mathcal F(X)$: 

For orbits 
$ \frak q, \widetilde{\frak q},\frak r, \widetilde{\frak r}\in \mathcal O^0(X)$
 we say that orbits
 $$
 \frak x \in \mathcal F(\frak q,\frak r), \quad  \quad 
 \widetilde {\frak x} \in \mathcal F(\widetilde {\frak q},  , \widetilde {\frak r}),
 $$
 are $\approx$-equivalent, if
 $$
 \frak q = \widetilde {\frak q}, \quad \quad  \frak r = \widetilde {\frak r},
 $$
 and if it holds for all choices of  
 $$
  \frak u \in \mathcal F(\bullet ,\frak q  ) ,  \quad  \quad \frak u \in \mathcal F(\frak r , \bullet ),
 $$
 that
 $$
\cnc \langle \frak u, \frak x, \frak v  \rangle \in \mathcal O(X), 
 $$
 if and only if
 $$
\cnc \langle \frak u,\widetilde{ \frak x}, \frak v  \rangle \in \mathcal O(X). 
 $$
 One checks that $\approx$ is an equivalence relation. By Lemma 2.1 it is invariantly attached to the subshift.
 
For orbits 
$ \frak q, \widetilde{\frak q},\frak r, \widetilde{\frak r}\in \mathcal O^0(X)$
 we say that orbits
 $$
 \frak x \in \mathcal F(\frak q,\frak r), \quad  \quad 
 \widetilde {\frak x} \in \mathcal F(\widetilde {\frak q},  , \widetilde {\frak r}),
 $$
 are $\cong$-equivalent, if
 $$
 \frak q \sim \widetilde {\frak q}, \quad \quad  \frak r \sim \widetilde {\frak r},
 $$
 and if it holds for all choices of  
 $$
 \widehat{\frak q}  \sim \frak q, \quad \quad  \widehat{\frak r}  \sim \frak r,
 $$
and of
  $$
 \frak u \in \mathcal F^0(\thinspace  \widehat{\frak q},\frak q),  \quad
 \widetilde {\frak u}
  \in \mathcal F^0(\thinspace  \widehat{\frak q},\widetilde {\frak q}), \quad
  \frak v \in \mathcal F^0(\thinspace  \widehat{\frak r},\frak r),  \quad
 \widetilde {\frak v}
  \in \mathcal F^0(\thinspace  \widehat{\frak r},\widetilde {\frak r}),
 $$
 that
 $$
 \cnc \thinspace \langle \frak u,\frak x,\frak v \rangle \approx 
 \cnc \thinspace\langle \widetilde {\frak u},\widetilde {\frak x},\widetilde {\frak v} \rangle.
$$
One checks that $\cong$-equivalence is an equivalence relation. By Proposition 3.1 
$\cong$-equivalence is invariantly attached to the subshift $X$. $\approx$-equivalence implies $\cong$-equivalence.

As a condition to be imposed on a subshift that satisfies $\bold H01-\bold H06$ we adopt the following  condition 

\smallskip
\noindent
$\bold H07$.
For $V \in \frak V(X) , \frak q, \frak r \in  V$, and for
$$
\frak x \in \mathcal F(\bullet, \frak q ),  \quad \frak x \in \mathcal F(\frak r, \bullet),
$$
and
$$
\frak u, \frak v \in \mathcal F^0(\frak q, \frak r ),
$$
one has that 
$$
\cnc \thinspace \langle  \thinspace {\frak x} ,  {\frak u}, \frak y \rangle  \in \mathcal O(X),
$$
if and only if
$$
\cnc \thinspace \langle  \thinspace {\frak x} ,  {\frak v}, \frak y \rangle \in \mathcal O(X).
$$

 \smallskip
 
We denote the set of $\cong$-equivalence classes by $\mathcal S(X)$. We give 
$\mathcal S(X)$ the structure of a semigroup (generally with zero). Let there be given 
$\cong$-equivalence classes $\mathcal D^{\langle - \rangle}$ and 
$\mathcal D^{\langle + \rangle}$ and  orbits
$$
\frak q^{\langle - \rangle} ,\frak r^{\langle - \rangle} , \frak q^{\langle + \rangle},
 \frak r ^{\langle + \rangle}   \in \mathcal O^0(X), 
$$
and 
$$
\frak x^{\langle - \rangle}  \in \mathcal F^0(\frak q^{\langle - \rangle} ,\frak r^{\langle - \rangle}  ) \cap \mathcal D^{\langle - \rangle},  \quad \quad
\frak x^{\langle + \rangle}  \in \mathcal F^0(\frak q^{\langle + \rangle} ,
\frak r^{\langle + \rangle}  ) \cap \mathcal D^{\langle + \rangle}.
$$
 If  
  $\frak r^{\langle - \rangle}  \not \sim \frak q^{\langle + \rangle} $ then we set the semigroup product of $\mathcal D^{\langle - \rangle}$ and
  $\mathcal D^{\langle + \rangle}$  equal to zero.  
  Otherwise we choose an orbit
  $\frak p \sim \frak r^{\langle - \rangle}$ and orbits
  $$
  \frak u^{\langle - \rangle}  \in \mathcal F^0(\frak r^{\langle - \rangle} , \frak p   ), 
  \quad      \quad \quad   
   \frak u^{\langle + \rangle}  \in \mathcal F^0( \frak p ,\frak q^{\langle + \rangle}  ),
  $$
   and set the semigroup product of $\mathcal D^{\langle - \rangle}$ and
  $\mathcal D^{\langle + \rangle}$  equal to the  $\cong$-equivalence class of an orbit
 for which we use the notation
  $$
  \frak z(\frak x^{\langle - \rangle} , \frak r^{\langle - \rangle}  , \frak u^{\langle - \rangle}  , 
  \frak p  , \frak u^{\langle + \rangle} , \frak q^{\langle + \rangle},\frak x^{\langle + \rangle}),
 $$ 
 and that we define as the orbit
 $$
  \cnc \langle \frak x^{\langle - \rangle},  \frak u^{\langle - \rangle} ,  \frak u^{\langle + \rangle} ,  \frak x^{\langle + \rangle}  \rangle.
 $$
 In this way the semigroup product will be well defined, provided that for choices of
 $$
\widetilde{\frak q}^{\langle - \rangle}  \sim \frak q^{\langle - \rangle},
\widetilde{\frak r}^{\langle - \rangle} \sim \frak r^{\langle - \rangle},
 \widetilde{\frak q}^{\langle + \rangle} \sim \frak q^{\langle + \rangle} ,
\widetilde{ \frak r }^{\langle + \rangle}  \sim { \frak r }^{\langle + \rangle}, 
 $$
 and  of
 \begin{align}
\widetilde{   \frak x       }^{\langle - \rangle}\cong \frak x^{\langle - \rangle} ,\quad \quad
\widetilde{   \frak x       }^{\langle + \rangle}\cong\frak x^{\langle + \rangle}, 
 \end{align}
  and choices of
  $$
 \widetilde{ \frak p} \sim \frak p,
  $$
  and of 
  $$
  \widetilde{ \frak u}^{\langle - \rangle}  \in \mathcal F^0(\widetilde{ \frak r }^{\langle - \rangle} ,  \widetilde{ \frak p} ), 
  \quad      \quad \quad   
   \widetilde{ \frak u}^{\langle + \rangle}  \in \mathcal F^0(  \widetilde{ \frak p} ,
   \widetilde{ \frak q }^{\langle + \rangle}  ),
  $$
  it holds that 
  $$
  \frak z(\frak x^{\langle - \rangle} , \frak r^{\langle - \rangle}  , \frak u^{\langle - \rangle}  , 
  \frak p  , \frak u^{\langle + \rangle} , \frak q^{\langle + \rangle},\frak x^{\langle + \rangle}) \cong
   \frak z(\widetilde{\frak x}^{\langle - \rangle} , \widetilde{\frak r}^{\langle - \rangle}  , 
   \widetilde{\frak u}^{\langle - \rangle}  , 
  \widetilde{\frak p } , \widetilde{\frak u}^{\langle + \rangle} , 
  \widetilde{\frak q}^{\langle + \rangle},\widetilde{\frak x}^{\langle + \rangle}),
 $$ 
 which is the same as 
 \begin{align}
  \cnc \langle \frak x^{\langle - \rangle},  \frak u^{\langle - \rangle} ,  \frak u^{\langle + \rangle} ,  \frak x^{\langle + \rangle}  \rangle \cong 
  \cnc \langle\widetilde{ \frak x}^{\langle - \rangle},  \widetilde{\frak u}^{\langle - \rangle} ,  \widetilde{\frak u}^{\langle + \rangle} ,  \widetilde{\frak x}^{\langle + \rangle}  \rangle.
  \end{align}
For the proof of (3.3) let
$$
\widehat{\frak q} \sim \frak q^{\langle - \rangle}, \quad \quad \widehat{\frak r} 
\sim \frak r^{\langle + \rangle},
$$
and let
 $$
 \frak v^{\langle - \rangle} \in \mathcal F^0(\widehat{\frak q},   \frak q^{\langle - \rangle}  ),  \quad \quad 
 \frak v^{\langle + \rangle} \in \mathcal F^0(\frak r^{\langle + \rangle} ,   \widehat{\frak r} ).
 $$
Set
$$
\widehat{\frak z}
 (\frak v^{\langle - \rangle},{ \frak q}^{\langle - \rangle} , \frak x^{\langle - \rangle}, \frak u^{\langle - \rangle}, \frak p,
 \frak u^{\langle + \rangle}, 
\frak x^{\langle + \rangle},{\frak r}^{\langle + \rangle},\frak v^{\langle + \rangle}) =
$$
$$
\cnc \langle \frak v^{\langle - \rangle},\frak x^{\langle - \rangle}, \frak u^{\langle - \rangle}, 
 \frak u^{\langle + \rangle}, 
\frak x^{\langle + \rangle},\frak v^{\langle + \rangle} \rangle.
$$
Also let
$$
 \widetilde{ \frak v}^{\langle - \rangle} \in \mathcal F^0(\widehat{\frak q},  
 \widetilde{ \frak q}^{\langle - \rangle}  ),  \quad \quad 
\widetilde{ \frak v}^{\langle + \rangle} \in \mathcal F^0(\widetilde{\frak r}^{\langle + \rangle} ,   \widehat{\frak r} ).
 $$
 To confirm (3.3) we prove that
 \begin{align*}
 &\widehat{\frak z}
 (\frak v^{\langle - \rangle},{ \frak q}^{\langle - \rangle} , \frak x^{\langle - \rangle}, \frak u^{\langle - \rangle}, \frak p,
 \frak u^{\langle + \rangle}, 
\frak x^{\langle + \rangle},{\frak r}^{\langle + \rangle},\frak v^{\langle + \rangle}) \approx
\\
&\widehat{\frak z}
 (\widetilde{\frak v}^{\langle - \rangle},\widetilde{ \frak q}^{\langle - \rangle} , 
\widetilde{ \frak x}^{\langle - \rangle},\widetilde{ \frak u}^{\langle - \rangle}, 
\widetilde{ \frak p},
\widetilde{  \frak u}^{\langle + \rangle}, 
\widetilde{ \frak x}^{\langle + \rangle},{\widetilde{ \frak r}}^{\langle + \rangle},
\widetilde{\frak v}^{\langle + \rangle}).
 \end{align*} 
 We choose auxiliary orbits
 $$
 \bar {\frak u}^{\langle - \rangle}  \in
  \mathcal F^0( \widetilde{ \frak r}^{\langle - \rangle}, \frak p  ), \quad \quad 
  \bar { \frak u}^{\langle + \rangle} \in 
  \mathcal F^0(\frak p , \widetilde{ \frak q}^{\langle + \rangle}  ).
 $$
 By (3.2) one has that
 $$
 \widehat{\frak z}
 (\frak v^{\langle - \rangle},{ \frak q}^{\langle - \rangle} , \frak x^{\langle - \rangle}, \frak u^{\langle - \rangle}, \frak p,
 \frak u^{\langle + \rangle}, 
\frak x^{\langle + \rangle},{\frak r}^{\langle + \rangle},\frak v^{\langle + \rangle}) \approx
 $$
 $$
  \widehat{\frak z}
 (\widetilde{\frak v}^{\langle - \rangle},{ \widetilde{\frak q}}^{\langle - \rangle} ,
 \widetilde{ \frak x}^{\langle - \rangle},
 \bar{ \frak u}^{\langle - \rangle}, \frak p,
 \bar{\frak u}^{\langle + \rangle}, 
\widetilde{\frak x}^{\langle + \rangle},{\widetilde{\frak r}}^{\langle + \rangle},
\widetilde{\frak v}^{\langle + \rangle}).
 $$
 and by $\bold H07$ one has that
 $$
  \widehat{\frak z}
 (\widetilde{\frak v}^{\langle - \rangle},{ \widetilde{\frak q}}^{\langle - \rangle} ,
 \widetilde{ \frak x}^{\langle - \rangle},
 \bar{ \frak u}^{\langle - \rangle}, \frak p,
 \bar{\frak u}^{\langle + \rangle}, 
\widetilde{\frak x}^{\langle + \rangle},{\widetilde{\frak r}}^{\langle + \rangle},
\widetilde{\frak v}^{\langle + \rangle}) \approx
 $$
 $$
 \widehat{\frak z}
 (\widetilde{\frak v}^{\langle - \rangle},\widetilde{ \frak q}^{\langle - \rangle} , 
\widehat{ \frak x}^{\langle - \rangle},\widetilde{ \frak u}^{\langle - \rangle}, 
\widetilde{ \frak p},
\widetilde{  \frak u}^{\langle + \rangle}, 
\widetilde{ \frak x}^{\langle + \rangle},{\widetilde{ \frak r}}^{\langle + \rangle},
\widetilde{\frak v}^{\langle + \rangle}).
 $$

One checks the associativity of this semigroup product.
 
For $V \in \frak V(X)$ and $\frak q, \frak r \in V$ the orbits in $\mathcal F^0(\frak q , \frak r  )$ are by $\bold H07$
$\cong$-equivalent. 
We denote their $\cong$ equivalence class by $\bold 1_{V}$.
\smallskip.
By $\bold H07$ one has for $U, W \in \frak V(X)$ and for $\frak q \in U, \frak r \in W$ and for 
$\frak y \in \mathcal F(\frak q , \frak r)$ that
$$
 \bold 1_{U}[\frak y]_{\cong} = [\frak y]_{\cong} = [\frak y]_{\cong}\bold 1_{W}.
$$
For $U, W \in \frak V(X), U \neq W,$ we have $UW= 0$ and one checks that  
$\bold 1_{V}, V \in \frak V(X)$, are idempotent.
We have the partition
$$
\mathcal S(X) = \bigcup_{U, W \in \frak V(X)}\bold 1_{U}\mathcal S(X) \bold 1_{W}.
$$
We have obtained local units $\bold 1_{V}, V \in \mathcal V(X),$ for $\mathcal S(X)$.

For topologically conjugate subshifts $X \subset \Sigma^\Bbb Z$ and 
$\bar X \subset \bar \Sigma^\Bbb Z$ that satisfy
$\bold H01$ - $\bold H07$ 
a topological conjugacy
$
\psi: X \to \bar X
$
carries the set $\mathcal S(X) $ with its equivalence relation $\cong$ into the set 
$\mathcal S(X) $ with its equivalence relation $\cong$. This means that the semigroup $\mathcal S(X)$ is invariantly associated to $X$. In particular
 $$
 \psi (\bold 1_{V}) = \bold 1_{\psi(V)}, \quad \quad V \in \frak V(X).
 $$

We denote by $\mathcal S^-(X)$($\mathcal S^+(X)$) the set of
$
f^- \in \mathcal S(X) \setminus \{0\} 
($
$
f^+\in \mathcal S(X) \setminus \{0\} 
$)
such that the following holds:
If $V \in \frak V(X)$ is such that $f^-=\bold 1_{V} f^-$
($f^-=f^+\bold 1_{V} $) then $0 \not \in \mathcal S(X)\bold 1_{V}f^-$
($0 \not \in f^+\bold 1_{V}S(X)$).

 We say  that $f^- \in \mathcal S^-(X) $($f^+ \in \mathcal S^+(X) $) is indecomposable if there do not exist $a^-, b^- \in \mathcal S^-(X)$($b^+, a^+ \in \mathcal S^+(X)$) such that
 $ f^-=a^-b^-$($ f^+=b^+a^+$). We denote the set of indecomposable elements of 
 $\mathcal S^-(X) $($\mathcal S^+(X) $) by $\mathcal E^-(X) $($\mathcal E^+(X) $).
 We set
 \begin{align*}
& \mathcal E^-(X)(U,W) = \bold 1_{U}\mathcal E^-(X)\bold 1_{W},
 \\
& \mathcal E^+(X)(U,W)= \bold 1_{U}\mathcal E^-(X)\bold 1_{W}, \quad 
\quad U, W \in \frak V(X).
 \end{align*}
 
 We have the partitions
 $$
 \mathcal E^-(X) = \bigcup_{U, W \in \frak V(X)} \mathcal E^-(X)(U,W),\quad
  \mathcal E^+(X) = \bigcup_{U, W \in \frak V(X)} \mathcal E^+(X)(U,W).
 $$
 
We state algebraic conditions $\bold H08$ -  $\bold H10$ on $\mathcal S(X)$ that allow to associate invariantly with the subshift $X\subset \Sigma^\Bbb Z$ an $\mathcal R$-graph.
The Condition $\bold H08$ and $\bold H09$ to be adopted come in two parts that are symmetric to each other.
Condition $\bold H10$ secures the existence of normal forms for the semigroup elements.
 
\smallskip
 
 \noindent
 $\bold H08$.
 The following holds for $U, W \in \frak V(X)$:
 
 \noindent
 $(-)$ For $e^- \in  \mathcal E^-(X)(U,W)$ there exists $e^+ \in  
 \mathcal E^+(X)(U,W)$ such that $e^-e^+ = \bold 1_{W}$.
 
  \noindent
 (+) For $e^+ \in  \mathcal E^+(X)(U,W)$ there exists $e^- \in  
 \mathcal E^-(X)(U,W)$ such that $e^-e^+ = \bold 1_{W}$.
  
  \smallskip
 
 \noindent
 $\bold H09$.
 The following holds for $U, W \in \frak V(X)$:
 
 \noindent
 $(-)$ For $f^-\in \mathcal S^-(X)(U,W)$ there exists uniquely $K \in \Bbb N$ , such that with the notation
 $$
 V_0 = U, \quad V_K = W
 $$
  there exist uniquely 
  $$
  V_k \in \frak V(X), 0 < k < K,
  $$
  and
 $$
 e^-_k\in\mathcal E^-(X)(V_{k-1},V_k), \quad 0 < k \leq K,
 $$ 
 such that
 $$
 f^- = \prod _{1\leq k \leq K}e^-_k.
 $$
 \noindent
 $(+)$
 For $f^+\in \mathcal S^+(X)(U,W)$ there exists uniquely $K \in \Bbb N$ , such that with the notation
 $$
 V_0 = W, \quad V_K = U
 $$
  there exist uniquely 
  $$
  V_k \in \frak V(X),  K > k > 0,
  $$
  and
 $$
 e^+_k\in\mathcal E^-(X)(V_{k+1},V_k), \quad K \geq  k > 0,
 $$ 
 such that
 $$
 f^+ = \prod _{K \geq  k > 0} e^+_k.
 $$
 \noindent

 \smallskip

 \noindent

\begin{flushleft}
$\bold H10.  \quad \quad 
\mathcal S(X) \setminus (\mathcal S^-(X) \cup \mathcal S^+(X)) = \bigcup_{V \in \frak V(x)} (\{\bold 1_{V}\} \cup
 \{\mathcal S^-(X) \bold 1_{V}\mathcal S^+(X)\})$.
 \end{flushleft}

For $U, W \in \frak V(X)$  we introduce a relation $\mathcal R(U,W)$ into
 $ \mathcal E^-(U, W) \times  \mathcal E^+(U, W)$ by $(e^-, e^+) \in \mathcal R(U,W)$
if and only if $e^-e^+ = \bold 1_{W}.$ We introduce a relation $\mathcal R(X)$ into
 $ \mathcal E^-(X) \times  \mathcal E^+(X)$ by
 $$
 \mathcal R(X) = \bigcup_{U,W \in \frak V(X)} \mathcal R(U,W).
 $$
 
 At this point we have obtained the $\mathcal R$-graph 
 $G_{\mathcal R(X)}( \frak V(X), \mathcal E^-(X)  \cup\mathcal E^+(X ))$.
 \begin{lemma}
 Let the  subshift $X\subset \Sigma^\Bbb Z$ satisfiy $\bold H01-\bold H10$.  Let $V \in  \frak V(X)$ and let $U$ be a predecessor vertex of $V$ in 
 $G( \frak V(X), \mathcal E^-(X) )$. Then
 
\noindent
$(-)$ For $e^-, \widetilde{e}^- \in \mathcal E^-(U,V)$ the equality of 
$
 \{ e^+ \in \mathcal E^+(X): (e^-,e^+) \in \mathcal R(X)\} 
 $
and
$
 \{ e^+ \in \mathcal E^+(X): (\widetilde{e}^-,e^+) \in \mathcal R(X)\} 
 $ 
implies the equality of $e^-$  and $ \widetilde{e}^-$.
 
\noindent
$(+)$ For $e^+, \widetilde{e}^+ \in \mathcal E^+(U,V)$ the equality of 
$
 \{ e^- \in \mathcal E^-(X): (e^-,e^+) \in \mathcal R(X)\} 
 $
and
$
 \{ e^- \in \mathcal E^-(X): (\widetilde{e}^-, \widetilde{e}^+) \in \mathcal R(X)\} 
 $ 
the equality of $e^+$  and $ \widetilde{e}^+$.
\end{lemma}
\begin{proof}
We prove $(-)$. Let $e^- \in \mathcal E^-(U,V),\widetilde{e}^- \in \mathcal E^-(U,V)$ and assume that
\begin{align}
 \{ e^+ \in \mathcal E^+(X): (e^-,e^+) \in \mathcal R(X)\} 
 =
 \{ e^+ \in \mathcal E^+(X): (\widetilde{e}^-,e^+) \in \mathcal R(X)\} ,
\end{align}
and
\begin{align}
e^- 
\neq \widetilde{e}^-.
\end{align}
Let $ \frak q,\widetilde{ \frak q} \in U,  \frak r,\widetilde{ \frak r} \in V$, and
\begin{align}
\frak x \in \mathcal F^-(\frak q , \frak r  ), \quad \widetilde {\frak x} 
\in \mathcal F^-(\widetilde{ \frak q} , \widetilde{ \frak r}  ),
\end{align}
such that
$$
[\frak x]_{\cong} = e, \quad \quad [\widetilde{\frak x}]_{\cong} = \widetilde{e}.
$$
By (3.6) 
$$
\frak x \cup \widetilde{\frak x} \subset A^-(X),
$$
and from (3.4) and $\bold H10$ it is seen that 
$$
[\frak x]_{\cong} = [\widetilde{\frak x}]_{\cong},
$$
contradicting (3.5).
\end{proof}

As a condition to be imposed on a subshift $X\subset \Sigma^\Bbb Z$ that satisfies 
 $\bold H01$ - $\bold H10$ we adopt the following Condition that comes in two parts 
$(-)$ and $(+)$ that are by $\bold H10$ 
and by Lemma 3.2
equivalent to each other.

\smallskip

 $\bold H11$. 
 
 $(-)$ The graph
 $G(\frak V(X),\mathcal E^-(X))$ is not a cycle.
 
 $(+)$ The graph
 $G(\frak V(X),\mathcal E^+(X))$ is not a cycle.
 
 \smallskip
 
 \begin{theorem}
 A subshift $X\subset \Sigma^\Bbb Z$ that satisfies $\bold H01-\bold H11$ is not synchronizing.
 \end{theorem}
 \begin{proof}
 
 By $\bold H11$ at least one of the following statements (a) and (b) holds:
 
(a) There exists a vertex $V \in  \frak V(X)$ with a single predecessor vertex $U$ in 
$G( \frak  V(X), \mathcal E^-(X))$, such that the set $ \mathcal E^-(U, V)$ contains at least two edges.

(b) There exists a vertex $V \in  \frak V(X)$ with at least two predecessor vertice $U$ and $W$ in $G( \frak  V(X), \mathcal E^-(X))$.

Assume (a). By Lemma 3.2 that there are edges 
$e^-\in \mathcal E^-(U, V), e^+\in \mathcal E^+(U, V),$ such that
\begin{align}
e^-e^+= 0.
\end{align}

Let $b \in \mathcal L(X)$. By  $\bold H07$  there exists $\frak p \in V$ such that the word $b$ appears on the points in $\frak p$. Let $X$ satisfy $\bold H03$ with parameter $M$ and let $m\geq M$ be such that $m  \times \card (\frak p)$ exceeds the length of $b$.  
Also let $\frak q \in U$ and let
$$
\frak x \in e^- \cap \mathcal F^-(\frak q, \frak p), \quad  \quad  \frak y \in e^+ 
\cap \mathcal F^+(\frak p, \frak q).
$$
By (3.9) 
$$
\cnc_m\langle \frak x, \frak y  \rangle \not \in \mathcal O(X).
$$
This means that the word $b$ cannot be synchronizing.
 
In case (b) choose edges 
$e^- \in \mathcal E^-(U,V), e^+ \in \mathcal E^+(W,V)$. Then
$$
e^-e^+ = 0.
$$
Proceed as in case (a).
\end{proof}

 \section{Hyposynchronization}
 
 We introduce notation that we use in this section and the next:
 For $p \in P^0(X)$ we denote by $F^-(p, \bullet)$($F^{(-)}(p, \bullet)$) the set of points in 
 $A^{-}(X)$($A^{(-)}(X)$)  that are left asymptotic to $p$. For $p, q \in P^0(X)$ we denote by $F^0(p, q)$($F^{+}(p, q)$) the set of points in 
 $A^{0}(X)$($A^{(+)}(X)$)  that are left asymptotic to $p$ and right asymptotic to $q$.
 
 \subsection{Hyposynchronizing shifts}
 
 In this subsection we consider a subshift 
 $X\subset\Sigma^{\Bbb Z}$ that satifies $\bold H02 - \bold H14$.
 We set
 \begin{align}
 \mathcal J_i(x) = \bigcap_{k\in \Bbb N}\{j \in (-\infty, i):
  \omega^{\langle + \rangle} (x_{(j-k,i ] })  = 
   \omega^{\langle + \rangle} (x_{(j , i ] })\}, \quad \quad i \in  \Bbb Z, x \in X,
 \end{align}
 and (in the case that $\mathcal J_i(x)$ is not empty) we set 
 \begin{align}
 J_i(x) =  \max \mathcal J_i(x), \quad \quad  i \in \Bbb Z, x \in X .
  \end{align}
 We set
 $$
 B(X) = \bigcap_{i\in \Bbb Z} \{x \in X : \mathcal J_i(x)  \not = \emptyset,
 x_{(i, \infty)} \in \omega_\infty^{\langle + \rangle} (x_{(J_i(x),i ] }  \}.
 $$
 By Lemma 2.1 the set $B(X)$ is invariantly attached to $X$.
 
 For $ i \in \Bbb Z, x^{( - \infty, i]} \in X_{( - \infty, i]}$ we give the notations 
 $\mathcal J_i(x^{( - \infty, i]})$ and $J_i(x^{( - \infty, i]})$ the meaning that is analogous to (4.1) and (4.2), and
 we set
  \begin{align*}
 B^{( - \infty, i]}(X) = \bigcap_{-\infty j \leq i}  &\{x^{( - \infty, i]} \in X_{( - \infty, i]} :
 \\
  &\mathcal J_j(x^{( - \infty, i]})  \not = \emptyset, \ \ 
   x^{( - \infty, i]}_{(j, i]} \in  ( \omega_\infty^{\langle + \rangle} 
 (x^{( - \infty, i]}_{(J_j(x^{( - \infty, i]}),j ] } ) )_{(j, i]}\}.
 \end{align*}
 One has that
 $$
 x_{(-\infty, i]} \in B^{( - \infty, i]}(X),         \quad \quad \quad x \in B(X), i \in \Bbb Z.
 $$
 We set
 \begin{align*}
& \omega_\infty^{\langle + \rangle}(x^{( - \infty, i]})  
 =  \omega_\infty^{\langle + \rangle}
 (x^{( - \infty, i]}_{(J_i(x), i]}),    \quad \quad 
  i \in \Bbb Z,x^{( - \infty, i]}\in B^{( - \infty, i]}(X).
\end{align*}

We adopt the following condition $\bold H12$ as a condition to be imposed on a subshift that satisfies $\bold H01 - \bold H11$

\smallskip

\noindent
$\bold H12$. For $p \in P^0(X)$ it holds that
$
F^{-}(p, \bullet) \subset B(X).
$

\smallskip

\noindent
Taking $\bold H12$ as the defining property 
we say that the subshift $X$ is hyposynchronizing if it satisfies $\bold H01 - \bold H12$.

 \smallskip
 
 \begin{proposition}
Let $\ell \in \Bbb N, i \in \Bbb Z$. Let 
 $$
x^{(-\infty ,i)} \in B^{(-\infty ,i)}(X),
 $$
 and
 \begin{align}
 a \in 
 \left(\omega_\infty^{\langle + \rangle}
 (x^{(-\infty ,i)})\right)_{[1, \ell]}.
  \end{align}
 Then
 \begin{align}
 (x^{(-\infty ,i)},a) \in B^{(-\infty ,i+ \ell)}(X),
  \end{align}
  and
  \begin{align}
  \omega_\infty^{\langle + \rangle}(x^{(-\infty ,i)}, a) =
  \left ( \omega_\infty^{\langle + \rangle}(x^{(-\infty ,i)}\right)_{(\ell, \infty)}
   \end{align}
  \end{proposition}
  \begin{proof}
  We show that
  \begin{align}
  J_i  (x^{(-\infty ,i)})    \in \mathcal J_{i + \ell}(x^{(-\infty ,i)},a),
  \end{align}
  which will confirm (4.2). For the proof of (4.6)
  we let $j < J_i(x^{(-\infty ,i)}) $ and
  $$
  u^{(i + \ell,\infty)} \in  \omega_\infty^{\langle + \rangle}
  ((x^{(-\infty ,i)}, a)_{[j, i+\ell]})
  $$
  and we provide a proof of
  \begin{align}
  u^{(i + \ell,\infty)} \in  \omega_\infty^{\langle + \rangle}
  ((x^{(-\infty ,i)}, a)_{[J_i (x^{(-\infty ,i)}, i+\ell]})
   \end{align}
   To set up this proof we let
   \begin{align}
   v^{(-\infty,J_i (x^{(-\infty ,i)}))} \in 
   \Gamma^{\langle - \rangle}_\infty (x^{(-\infty ,i)}, a)_{[J_i (x^{(-\infty ,i)}, i+\ell]}  ).
   \end{align}
 The task is then to show that
 \begin{align}
  v^{(-\infty,J_i ( v^{(-\infty,J_i (x^{(-\infty ,i)}))}))}x^{\langle - \rangle}_{[J_i(x^{(-\infty ,i)}), 0]}au^{(i + \ell,\infty)} \in X.
 \end{align}
 By (4.3)
 \begin{align}
 a \in  \left (  \omega_\infty^{\langle + \rangle}( v^{(-\infty,J_i (x^{(-\infty ,i)}))}_{[j, 0]})\right)_{[1, \ell]}.
 \end{align}
 By (4.7) and (4.10)
 $$
 au^{(i + \ell,\infty)} \in  \omega_\infty^{\langle + \rangle}(x^{(-\infty ,i)}_{[j, 0]}).
 $$
 Therefore
 $$
 au^{(i + \ell,\infty)} \in 
  \omega_\infty^{\langle + \rangle}(x^{(-\infty ,i)}_{[J_i( x^{(-\infty ,i)}), 0]}),
 $$
 which means that (4.9) holds.
 \end{proof}
  \begin{corollary}
 Let  
 $$
x^{(-\infty ,i)} \in B^{(-\infty,i)}(X),
 $$
 and let
 $$
 x^{( i,\infty)} \in \omega_\infty^{\langle + \rangle}(x^{(-\infty ,i)]} ).
 $$
 Then
 $$
 (x^{(-\infty ,i)},x^{( i,\infty)}  ) \in B(X).
 $$
 \end{corollary}
 \begin{proof}
 Apply the proposition to
 $$
 a = x^{( i,\infty)} _{(i, i+ \ell]},  \quad \quad \ell \in \Bbb N. \qed
 $$
  \renewcommand{\qedsymbol}{}
 \end{proof}
 
 As a condition to be imposed on a hypersynchronizing subshift we adopt the following
 condition $\bold H13$:
 
 $\bold H13$
 There exists $M \in \Bbb N$ such that the following statement holds:
For $V \in \frak V(X)$, and $p, q \in V$, let 
$x \in F^0(p, q)$  such that
$$
I^+(x) < M \per (q).
$$
Then
 $$
  \omega_\infty^{\langle + \rangle} ( {x}_{(-\infty, 0 )} )= \vrt (q).
  $$

\smallskip

\noindent

 For $\frak p \in \mathcal O^0(X)$ we can define by  $\bold H13$ a vertex set 
 $\frak V_\frak p(X)$ by 
 $$
 \frak V_\frak p(X) = \bigcup_{p \in \frak p, j \in \Bbb Z} 
 \{\omega^{\langle + \rangle}( y_{( - \infty, j]} ) : y \in F^{(-)} (p, \bullet)\}.
 $$

 \begin{lemma}
 Let the hyposynchronizing subshift $X$ satisfy  $\bold H13$ with parameter $M$.
 Let $V \in \frak V(X),$ and let $\frak p, \frak q \in V$. Then
 $\frak V^{hypo}_{\frak p}(X)$ equals $\frak V^{hypo}_{\frak q}(X)$.
  \end{lemma}
  \begin{proof}
  Let $q \in \frak q$ and let $y \in F^{(-)}(q, \bullet)$. Also let $j\in \Bbb Z$.
  We prove that
  $$
  \omega^{\langle + \rangle}( y_{( - \infty, j]} \in \frak V_\frak p(X).
  $$
  Without loss of generality we can assume that
  $$
  j > I^-(y)  - \per(q).
  $$
  There exist $p \in\frak p$ and $x \in F^0(p, q)$ such that
$$ 
I^+(x)    < I^-(y)  - (M+1)\per(q).
$$ 
Set
$$
z = (x_{( - \infty, I^-(y) - (M+1)\per(q))}, y_{ [I^-(y) - (M+1)\per(q), \infty)}).
$$
One has that $z\in  F^{(-)}(p, \bullet)$, and by $\bold H13$ one has that
$$
\omega^{\langle + \rangle}(z_{( - \infty, I^-(y) - \per(q)]}) 
= \omega^{\langle + \rangle}(y_{( \infty,I^-(y) - \per(q) )}).
$$
Therefore also
$$
\omega^{\langle + \rangle}(z_{( - \infty, j]}  ) =\omega^{\langle + \rangle}(y_{( - \infty, j]}  ) \qed
$$
 \renewcommand{\qedsymbol}{}
 \end{proof}
  As notation for the common value of
   $\frak V^{hypo}_\frak p(X), \frak p \in V,$ we use the symbol $\frak V^{hypo}_V(X)$.
 We set
 $$
 \frak V^{hypo}(X) = \bigcup_{V \in\frak V^{hypo}(X)} \frak V^{hypo}_V.
 $$
 
 As a conditon to be imposed on a hyposynchronizing subshift that satisfies $\bold H12$ and 
 $\bold H13$ we adopt
 
 $\bold H14$.There exists $M \in \Bbb N$ such that the following statement holds: For $V \in \frak V(X)$, let  there be given a finite path
\begin{align}
f^- = ( f^-_k)_{1\leq  k \leq K}
\end{align}
in $G_V^{hypo}(X)$
such that $\src (f^- ) = V$ together with a neutral periodic orbit
$\frak p \in \trg(f^-)$
and a neutral periodic point 
 \begin{align}
p \in\frak p.
\end{align} 
Let
$$
\frak q_k \in \src(f^-_k),\quad \quad 1 \leq k \leq K,
$$
and let
$$
\frak y_K \in \mathcal F^-(\frak q_K , \frak p )\cap f^-_K,
$$
$$
\frak y_k \in \mathcal F^-(\frak q_k , \frak q_{k+1}  ) \cap f^-_k, \quad \quad 1 \leq k < K .
$$
Let 
$$
x \in \cnc\langle \frak y_k \rangle_{1 \leq k \leq K} \cap F(\bullet,p)
$$
such that
$$
I^+(x) < M.
$$
Then the vertex
$\omega^{\langle + \rangle}_\infty ( x_{(- \infty, 0]} )$
of $G_{hypo}(X))$ is uniquely determined by the choices (4.13) of $f^- $ and (4.14) of $p$.

We denote this vertex that is produced by $\bold H14$  by $V^{hypo}(f^- , p) $.

  One one can read off from the proposition and corollary the data of a 
 $\Sigma$-labelled one-right resolving directed graph
 $
 G^{hypo}_V(X)
 $
 with vertex set $\frak V^{hypo}_V(X)$:
The set $\Sigma(V^{hypo} )$ of symbols that are accepted by the vertex 
$V^{hypo} \in  \frak V^{hypo}_V(X)$ is given by
the set of symbols $\sigma \in \Sigma$ that appear as initial symbols of points in
$V^{hypo}$. The target vertex of the edge that leaves $V^{hypo}$ and that carries the symbol $\sigma \in \Sigma(V^{hypo})$ is the vertex that contains the points that are obtaines by removing from points in $V^{hypo}$ that have $\sigma$ as initial symbol this initial symbol $\sigma$. 

In the same way one has a $\Sigma$-labelled one-right resolving directed graph
 $
 G^{hypo}(X)
 $
 with vertex set $\frak V^{hypo}(X)$.

 It follows from Lemma 4.1 that
 the graph $G^{hypo}_V(X)$ 
 is strongly connected, and by $\bold H08$ it presents $X$.

  \subsection{Edge shifts}

 \begin{lemma}
 Let$X \subset \Sigma^\Bbb Z$ and 
 $\bar X \subset \bar\Sigma^\Bbb Z$ hyposynchronizing shifts. Let 
 $$ 
 \varphi: X \to \bar X
 $$ 
 be a topological conjugacy that is implemented by a block map $\Phi$ with coding window $[-L, L]$. Let
 $$
 \textcircled{\tiny $\bar V$} = \varphi( \textcircled{\tiny V}),
 $$ 
 and for points $x \in X$ use the notation
 $$
 \bar x = \varphi(x).
 $$
Then one has  a  topological  conjugacy $\varphi_E$ of
 $
 E(G^{hypo}_{\textcircled{\tiny V}} (X))$  onto  
 $E(G^{hypo}_{ \textcircled{\tiny $\bar V$}}(X)), 
 $
 that is implemented by a block map $\Phi_E$ with memory $-3L$ and anticipation $L$  that is given by
$$
\Phi_E((V_i, x_i)_{-3L \leq i \leq L}) =
$$ 
 \begin{align*}
 (\{\bar y^{\langle + \rangle}_{[0, \infty)}: \bar y^{\langle + \rangle} \in \Phi 
( \{(x_{[-3L, -2L)}, x^{\langle + \rangle}) : x^{\langle + \rangle} \in V_{-L}  \}),   
 \bar y^{\langle + \rangle}_{[-2L, 0)} = \bar x_{[-2L, 0)} \}, 
\\
 \Phi (x_{[-L, L]}) ). 
  \end{align*}

 \end{lemma}
 \begin{proof}
 Setting hypothetically
 $$
 \varphi_E((V_i, x_i)_{-\infty < i < \infty}) = (\bar V_i,\bar x_i)_{-\infty < i < \infty},
 \quad\quad\quad (V_i, x_i)_{-\infty < i < \infty} \in  E(G^{hypo}_{\textcircled{\tiny V}} (X)),
 $$ 
 one arrives at this expression for $\bar V_0^{hypo}$ by collecting the points in 
 $V_{-L}^{hypo} $ whose image under $\Phi$ contributes to  $\bar V_0^{hypo}$.
 \end{proof}
 
 The lemma shows that the system of edge shifts 
 $E(G_V^{hypo}(X))$, and in particular also the edge shift
 $
 E(G^{hypo}(X))
 $,
  are invariantly attached to $X$. Compare the construction of extensions in \cite{Kr2, Kr3, Kr7}.

 \subsection{Cycles}
 
 Consider a hyposynchronizing subshift $X \subset \Sigma^\Bbb Z$ that satisfies 
 $\bold H12$ - $\bold H15$. Let $V \in \frak V(X)$. We say that a finite path 
 \begin{align}
 c = 
 (V^{hypo}
 _i, x_i)
 _{0 \leq i < \ell}
 \end{align}
 in the graph $G_{V}^{hypo}(X)$ is a cycle (more precisely, a cycle of length $\ell(c) = \ell$ at the vertex $V^{hypo}_0$), if
  \begin{align}
 \tau_{x_{i-1}}(V^{hypo}_{i-1}) =V^{hypo}_{0}.
  \end{align}
  
  A cycle $c$ as in (4.9) and (4.10)  determines a periodic point 
  $$
  p^E(c) = (V^{hypo}_j, x_j)_{j \in \Bbb Z}
  $$
  in $E(G^{hypo}_V)$
  of period $\ell(c)$ by
  $$
  (V^{hypo}_{k\ell(c)+i}, x_{k\ell (c)+i})_{0 \leq i < \ell(c)} = (V_i, x_i) _{0 \leq i < \ell(c)} \quad \quad 
  k \in \Bbb Z.
  $$
  The label sequence of this  periodic point we denote by
  $$
  p(c) = (x_j)_{j \in \Bbb Z}.
  $$
  
  We adopt as a condition to be imposed on a hypersynchronizing subshift that saisfies 
  $\bold H12 - \bold H14$ the following Condition
  
  \smallskip

  \noindent
 $\bold H15$.
For a cycle $c$ in $G^{hypo}_V, V \in \frak V(X)$ one has that $p(c)$ is a neutral periodic point.

\smallskip  

A neutral periodic point $p$ determines a cycle $c(p)$ of length $\per (p)$ by
$$
c(p) = (\omega^{\langle + \rangle} ( p_{(- \infty, i}), p_i)_{0\leq i < \per (p)}.
$$
One has that
$$
p(c(p)) = p, \quad \quad p\in P^0(X).
$$
We say that a cycle $c$ is zero-height if $c = c(p(c))$.

We adopt as a condition to be imposed on a hypersynchronizing subshift that satisfies
$\bold H12$ - $\bold H15$ the following condition $\bold H16$

$\bold H16$. Let $c$ be a cycle, let $p$ be the neutral periodic point of $c$. Then either
$$
c = (\omega (p_{( -\infty, i]} ), p_i )_{0 \leq i < \per(p)},
$$
or there exists a finite path $f^-$ in  $G(\frak V(X), \mathcal E^-(X)$ and $\frak p,  \in \trg( f^-), p\in \frak p$ such that
$$
c = ( V^{hypo}(f^-, p), p_i)_{0 \leq i < \per(p)}. 
$$ 
By Lemma 3.2 the path $f^-$ of $\bold H16$ is unique. The height of the cycle $c$ is defined as the length of $f^-$.

We adopt as a condition to be imposed on a hypersynchronizing subshift that saisfies 
  $\bold H12 - \bold H16$ the following Condition
\noindent

$\bold H17$.
For $p \in P^0(X)$, for $x \in   F^{(-)}(p, \bullet)$, and for $i \in \Bbb Z$ there exists 
$y \in F^0(p, p)$ such that
$
x_{(-\infty, i]} = y_{(-\infty, i]}.
$

\begin{lemma}
Every vertex of 
$G_{\textcircled{\tiny V}}^{hypo}(X) $
is on a cycle of zero-height.
\end{lemma}
\begin{proof}
Let $V^{hypo} \in \frak V_{\textcircled{\tiny V}}^{hypo}(X)$, and let 
$p \in \textcircled{\tiny V}$ and $x \in F^{(-)}(p, \bullet)$ such that
$$
V^{hypo}= \omega^{\langle + \rangle} ( x_{(- \infty, 0]}).
$$
By $\bold H16$ an $\bold H17$
there exists $y \in F^{0}(p, p)$ such that
$$
x_{(- \infty, 0]}  =  y_{(- \infty, 0]}  ,
$$
and such that for some $M \in \Bbb N$,
$$
\omega^{\langle + \rangle} ( p_{(- \infty, 0)}) = 
\omega^{\langle + \rangle} ( y_{(- \infty, I^+(y) + M \per(p))}).
$$
A point $z\in F^-(p, \bullet)$ is given by
$$
z = ( y_{(- \infty, I^+(y) + M \per(p))},   x_{(I^-(x), \infty)}).
$$
The point $z$ yields a cycle $c$ of zero-height of length
$$
\ell = I^+(x)- I^-(x)+ M \per(p) 
$$
by
$$
c = \omega^{\langle + \rangle} _\infty( z_{(- \infty, i)}, z_i)_{0 \leq i < \ell}.
$$
One has that
$$
V^{hypo}=\omega^{\langle + \rangle}_\infty ( p(c)_{(- \infty, 0)}). \qed
$$
 \renewcommand{\qedsymbol}{}
  \end{proof}

\noindent

\smallskip
\noindent

\bigskip

\smallskip
\noindent

As a condition to be imposed on a hyposynchronizing subshift that satisfies 
$\bold H12$ - $\bold H17$ we adopt the following Condition

$\bold H18$. For every vertex $V^{hypo}$ of $G_{\textcircled{\tiny V}}^{hypo}(X)$ there is a bound on the heights of cycles at $V^{hypo}$.

\begin{lemma}
Let $X$  be a hyposynchronizing subshift.

The following statements (a) and (b) are equivalent:

\noindent
(a) There exists an edge $e$ 
and cycles 
$$
c^{(k)}= (e^{(k)}_l)_{0 \leq l_\circ < \ell(c^{(k)} )}  ,   \quad \quad k \in \Bbb N,
$$ 
in
$G_{\textcircled{\tiny V}}^{hypo}(X) $ 
such that
\begin{align*}
e = e^{(k)}_0, \quad \quad \quad \eta(  c^{(k+1)} ) > \eta(c^{(k)} ) ,  \quad  k \in \Bbb N. 
\end{align*}
(b) For $l \in \Bbb N$ there exists in $G_{\textcircled{\tiny V}}^{hypo}(X) $ a path 
$a  = (e_{l_\circ}^{(k)})_{0 \leq l_\circ < l}$ 
 and cycles
$$
c^{(k)}= (e^{(k)}_{l_\circ})_{0 \leq  l_\circ< \ell(c^{(k)} )}  ,\quad \quad  k \in \Bbb N, 
$$ 
such that
$$
\ell(c^{(k)})> l, \quad \quad a = e^{(k)}_{0 \leq {l_\circ} < l}, \quad \quad \eta(c^{(k+1)} ) > \eta(c^{(k)}),  \quad   k \in \Bbb N.
$$
\end{lemma}
\begin{proof}
Assume (a). Let $l \in \Bbb N$ and set
$$
k_\circ = \max\{k \in \Bbb N : \ell(c^{(k)})  \leq l\}.
$$
There is a path that appears infinitely often in the sequence
$$
(e^{(k)}_{l_\circ})_{0 \leq l_\circ <  l },\quad \quad k > k_\circ.
$$
\end{proof}

\begin{proposition}
$\bold H18$  is an invariant of topological conjugacy
\end{proposition}
\begin{proof}
Statement (a) of Lemma 4.4 is the negation of $\bold H18$. In view of  Lemma 4.4 we prove that its statement (b) is invariant under topological conjugacy.

 Let there be given subshifts $X \subset \Sigma^\Bbb Z$, $\bar{X} 
 \subset \bar\Sigma^\Bbb Z$, and a topological conjugacy
 $$
 \psi: \bar{X} \to X
 $$
 that is implemented by a block map
 $$
 \Psi: \bar{X}_{[-L, L]} \to \Sigma.
 $$
 Let 
 $$
 \psi: E(G_{ \textcircled{\tiny $\bar V$}}^{hypo}(\bar X)) \to E(G_{\textcircled{\tiny V}}^{hypo}(X))
 $$
 be the conjugacy
 implemented by the block map
 $$
 \Psi^E: E(G^{hypo}_{ \textcircled{\tiny $\bar V$}}(\bar X))_{[-3L, L]} \to E(G_{\textcircled{\tiny V}}^{hypo}( X)),
 $$
 that is given by  Lemma 4.2.
We assume that $\bar X$ satisfies  statement (b) of Lemma 4.4 and we prove that $X$ also satisfies  
 statement (b) of Lemma 4.4:
 Let $l \in \Bbb N$, let 
 $ \bar a  = ( \bar a_{l_\circ}^{(k)})_{0 \leq l_\circ  < l+ 2L}$ be a path in 
 $G^{hypo}_{ \textcircled{\tiny $\bar V$}}(\bar X)$ and let 
$$
 \bar c^{(k)}= (\bar c^{(k)}_l)_{0 < l < \ell(c^{(k)} )}  ,\quad \quad  k \in \Bbb N, 
$$ 
be cycles in $G_{ \textcircled{\tiny $\bar V$}}^{hypo}(\bar X)$
 such that
 $
 \ell( \bar c^{(k)})> l + 2L,
 $
 and such that
\begin{align}
 \bar a =  (\bar c^{(k)})_{0 \leq {l_\circ} < l + 2L}, \quad \quad \eta(   \bar c^{(k+1)} ) > \eta( \bar c^{(k)} ),  \quad 
  k \in \Bbb N.
\end{align}
 By  Lemma 4.2  
 the  path $a = \Psi( \bar a)$ and the cycles $(c^{k})$ that are given by
 $$
 p^E(c^{(k)}) = S_\Sigma^L(\psi (p^E( \bar c^{(k)} ) ).  \quad   k \in \Bbb N,
 $$
 satisfy statement (b) of Lemma 4.4.
 \end{proof}
 
For a vertex $V^{hypo}\in \frak V_{\textcircled{\tiny V}}^{hypo}$ we define its height $\eta(V^{hypo})$
 in $G_{\textcircled{\tiny V}}^{hypo}(X) $ as the maximal height of cycles in $G_{\textcircled{\tiny V}}^{hypo}(X) $ at $V^{hypo}$. 
In view of $\bold H16$
the same notion of a height in $G_{\textcircled{\tiny V}}^{hypo}(X) $ of a vertex 
$V^{hypo}\in \frak V_{\textcircled{\tiny V}}^{hypo}$
can be introduced in the following way:

The height of a vertex $V^{hypo}\in \frak V_{\textcircled{\tiny V}}^{hypo}$
in $G_{\textcircled{\tiny V}}^{hypo}(X) $ equals the maximal length of a path $h^+$ in 
$G_{\textcircled{\tiny V}}^{hypo}(X) $ for which there are 
$$
\frak p \in \src(h^+), p \in \frak p,
$$ 
$$
\frak q  \in \textcircled{\tiny V}, q \in \frak q,
$$
and an orbit 
$$
\frak x \in      h^+ \cap \mathcal F^+(\frak p , \frak q  ),
$$
that contains a point
$$
x \in \frak x \cap F^+(p, q), 
$$
such that $I^-(x) < 0$ and such that
$$
x_{(I^-(x), \infty)} \in V^{hypo}.
$$

The vertex $V^{hypo}\in \frak V_{\textcircled{\tiny V}}^{hypo}$ is of  zero height in $G_{\textcircled{\tiny V}}^{hypo}(X) $
if and only if  
 no such path exists.

 \section{Restricted complexity}

 \subsection{A family of hyposynchronizing shifts with restricted complexity }
  
  We introduce a class of 
 hyposynchronizing subshifts that have restricted complexity. The starting point is an 
 $\mathcal R$-graph with a single vertex $\textcircled{\tiny V}$, and an edge set $\mathcal E^- \cup \mathcal E^+$ that carries a complete binary relation $\mathcal R$ that satisfies the following conditions $A(-)$ and $A(+)$ that are taken from Lemma 3.4:

 \smallskip
 
 $A(-)$:
 For $e^-, \widetilde {e}^- \in \mathcal E^-$ the equality of
 $ \mathcal E^+(e^-)$ and $ \mathcal E^+(\widetilde {e}^-)$
 implies the equality of
 $
 e^-
 $
 and
 $
   \widetilde {e}^-.
 $
 
$A(+)$: 
For $e^+, \widetilde {e}^+ \in \mathcal E^+$ the equality of
 $
 \mathcal E^-(e^+)
$
and
$
\mathcal E^-(\widetilde {e}^+)
 $
 implies the equality of
 $
 e^+
 $
 and
 $
  \widetilde {e}^+.
 $

 \smallskip
 
Let there be given a mapping
$$
\mu^-: \mathcal E^- \times \Bbb N \to    \Bbb Z_+,
$$
such that for $e^-\in \mathcal E^-$
$$
0 < \sum_{\lambda \in \Bbb N}\mu^-( e^-,\lambda) < \infty
$$
and a mapping 
$$
\mu^+: \mathcal E^- \to    \Bbb N.
$$
 Set for  $e^-\in \mathcal E^-$
\begin{align*}
& \Sigma^-( e^-) = 
 \{ (e^-,\lambda,
 \mu^- _{\circ}):
 \lambda \in \Bbb N, 1 \leq  \mu^- _{\circ} \leq    \mu^-( e^-,\lambda)\},
\\
& \Sigma^+( e^-) = 
 \{ (e^+,
 \mu^+ _{\circ}): e^+ \in \mathcal E^+(e^-),  1 \leq  \mu^+ _{\circ} \leq    \mu^+( e^+)\}, \quad \quad  \quad .
 \end{align*}

 Denote by
 $$   
  \mathcal A_{\mu^-, \mu^+}(G_{\mathcal R}(\{\textcircled{\tiny V} \}, \mathcal E^- \cup \mathcal E^+))
 $$
 the pushdown automaton with input alphabet
  $$
 \{ \Sigma^-(e^-)\cup \Sigma^+(e^-):       e^-\in    \mathcal E^-\}
 $$
and stack alphabet  $\mathcal E^-$,   
and with a list of transitions that is reduced to the following list:
If the stack of the automaton is empty then upon accepting the symbol
$$
 (e^-,\lambda,
 \mu^- _{\circ})  \in \mathcal E^-(e^-), \quad (e^-\in    \mathcal E^-),
$$
 as input, the automaton puts the word $(e^-)^\lambda$ into its stack.
 
 If the stack of the automaton contains a word
 $
 b
 $
 in the symbols of $ \mathcal E^-$,
 then upon accepting the symbol
 $$
 (e^-,\lambda,
 \mu^- _{\circ})  \in \mathcal E^-(e^-), \quad (e^-\in    \mathcal E^-),
$$
 as input, the automaton changes its stack content to $b(e^-)^\lambda $.

 For a word $b$
 in the symbols of $ \mathcal E^-$, and for $e^-\in    \mathcal E^-$, if the stack 
 of the automaton contains the word $be^-$, then upon accepting the symbol 
 $$
 (e^+,\mu^+ _{\circ}) \in \Sigma^+(e^-)
 $$ 
 as input the automaton changes its stack content to $b$.

 If the stack of the automaton contains the symbol $e^- \in  \mathcal E^-$ then, upon accepting the symbol 
 $$
  (e^+,\mu^+ _{\circ}) \in\Sigma^+( e^-)
  $$
  as input, the automaton empties its stack.
 
 The subshift that is presented by the automaton 
 $\mathcal A_{\mu^-, \mu^+} (G_{\mathcal R}(\{ \textcircled{\tiny V}\}, \mathcal E^- \cup \mathcal E^+))$ we denote by
 $Y_{\mu^-, \mu^+} (G_{\mathcal R}(\{ \textcircled{\tiny V}\}, \mathcal E^- \cup\mathcal E^+))$, and we denote the family that contains the subshifts $Y_{\mu^-, \mu^+} (G_{\mathcal R}(\{ \textcircled{\tiny V}\}, \mathcal E^- \cup\mathcal E^+))$ by $\mathcal Y$. 
 By construction the the subshifts 
 $Y_{\mu^-, \mu^+}( G_{\mathcal R}(\{ \textcircled{\tiny V}\}, \mathcal E^- \cup\mathcal E^+))$ have restricted complexity and
 $$
 G^{hypo}(Y_{\mu^-, \mu^+}
 ( G_{\mathcal R}(\{ \textcircled{\tiny V}\}, \mathcal E^- \cup\mathcal E^+)))=
 \mathcal A_{\mu^-, \mu^+} (G_{\mathcal R}(\{ \textcircled{\tiny V}\}, \mathcal E^- \cup \mathcal E^+)).
 $$

 It follows from Statements $A(-)$ and $A(+)$ for $\mathcal R$-graphs  
 $G_{\mathcal R}(\{\textcircled{\tiny V}\}, \mathcal E^- \cup  \mathcal E^+) )$  
with mappings  $\mu^-$,   $\mu ^+ $   
 and 
 $G_{\widetilde{\mathcal R}}
 (\{\widetilde{\textcircled{\tiny V}}\},
{\widetilde{ \mathcal E}}^- \cup {\widetilde{ \mathcal E}}^+) ), $ 
with mappings $\widetilde{ \mu}^-, \widetilde{\mu }^+ $,
 that the shifts 
 $Y_{\mu^-, \mu^+}(G_{\mathcal R}(\{ \textcircled{\tiny V}\}, \mathcal E^- \cup \mathcal E^+))$ 
 and $Y_{\widetilde{\mu}^-, \widetilde{\mu}^+}( G_{\widetilde{\mathcal R}}
 (\{\widetilde{\textcircled{\tiny V}}\},
{\widetilde{ \mathcal E}}^- \cup {\widetilde{ \mathcal E}}^+) )$
 are topologically conjugate if and only if there exists an isomorphism of partitions
$$
\chi:  {\widetilde  {\mathcal E}}^-\cup 
{\widetilde {\mathcal E}}^+
 \to \mathcal E^- \cup \mathcal E^+
$$
such that
 $$
 \widetilde  {\mu} ^-=\mu^- \circ (\chi \times \id ),  \qquad  \widetilde  {\mu}^+=\mu^+ \circ \chi. 
 $$
 
 In the case that
 $$
 \card (\mathcal E^-) = \card (\mathcal E^+),
 $$
 and
 $$
 \card ( \mathcal E^+(e^-)) = 1, \qquad (e^- \in \mathcal E^- ),
 $$
 $$
 \card ( \mathcal E^-(e^+)) = 1, \qquad (e^+ \in \mathcal E^+ ).
 $$
 the subshifts
 $  
  Y_{\mu^-, \mu^+}
  ( G_{\mathcal R}(\{ \textcircled{\tiny V}\}, \mathcal E^- \cup \mathcal E^+))
 $  
are  contained in the family of subshifts that was introduced in \cite{BH2} (see also \cite {D, LY, JM}).

\subsection{A literal conjugacy}

We  consider a hyposynchronizing subshift $X\subset \Sigma^\Bbb Z$ that satisfies $\bold H13$ - $\bold H18$. We adopt as a condition to be imposed on a hyposynchronizing shift that satisfies $\bold H13$ - $\bold H18$ the following Condition 
  
$\bold H19$.  $\card(\frak V(X)) =1$.

For a hyposynchronizing shift $X$ that satisfies $\bold H13$ - $\bold H19$ denote the set of fixed points of $X$ by $P_1(X)$, and
 set
 $$
 P_1^-(X) =P_1(X) \cap A^-(X), \quad P_1^+(X) = P_1(X) \cap A^+(X).
 $$
 We denote the set of symbols that are carried by the points in $P_1^-(X)$
 ($P_1^+(X)$) by $\Xi^-(X)$($\Xi^+(X)$).
 Set 
 $$
 P^0(V^{hypo}) = \{p \in P^0(X) :  \tau_{p_{[0, \per(p))}} (V^{hypo}) = V^{hypo} \},
 \quad\quad\quad         V^{hypo} \in\frak  V^{hypo}(X),
 $$
 $$
 P^0(\xi^-) = \bigcup_{m \in \Bbb N}\{p \in P^0(X): \xi^- \in \omega^{\langle + \rangle} 
 (p_{[- m \times \per (p))}) \}, \quad\quad\quad    \xi^- \in \Xi^-(X),
 $$
 $$
 P^0(\xi^+)= \bigcup_{m \in \Bbb N}\{p \in P^0(X): \xi^+ \in \omega^{\langle - \rangle} 
 (p_{[- m \times \per (p))}) \}, \quad\quad\quad    \xi^+\in \Xi^-(X),
 $$
  \begin{align*}
 P_\bullet^0(X) =\left (\bigcap_{V^{hypo} \in \frak V^{hypo}(X)}   P^0(V^{hypo} )\right )  &\cap
 \\ 
 \left (\bigcap_{\xi^- \in \Xi^-(X)}   P^0(\xi^-)\right ) &\cap
  \left (\bigcap_{\xi^+ \in \Xi^-(X)}   P^0(\xi^+)\right ). 
 \end{align*}
 
 The condition 
  \begin{align*}
 P_\bullet^0(X) \neq \emptyset.
 \end{align*}
 is not invariant under topological conjugacy (See subsection 6.2).
 \begin{lemma}
 Let $X \subset \Sigma^\Bbb Z$ be a hyposynchronizing subshift that satisfies 
 $\bold H13$ - $\bold H19$, such that
 \begin{align*}
 P_\bullet^0(X) \neq \emptyset.
 \end{align*}
 Let $q, r \in P_\bullet^0(X)$. Then
 
$(-)$ For $\xi^- \in \Xi^-(X)$,
$$
q_{(-\infty, 0)}\xi^-r_{[0, \infty)}  \in F^-(q,r).
$$

$(+)$ For $\xi^- \in \Xi^+(X)$,
$$
r_{(-\infty, 0)}\xi^+q_{[0, \infty)}  \in F^+(r,q).
$$
\end{lemma}
\begin{proof}
 $(-)$
 Let $c \in \Gamma^{\langle- \rangle}(\xi^-)$. There is a vertex in $G_{\textcircled{\tiny V}}^{hypo}(X)$ that is the target vertex of a path in $G_{\textcircled{\tiny V}}^{hypo}(X)$ with label sequence $c\xi^-$. It follows that the sequence
 $c\sigma r_{[1, \infty)}$ is the label sequence of a right infinite path in $G_{\textcircled{\tiny V}}^{hypo}(X)$.
 This  proves that
 \begin{align}
 r_{[0, \infty)}\in  \omega_\infty^{\langle+ \rangle} (\xi^-).
 \end{align}
 
One has by the choice of $p$
an $m \in \Bbb N$ such that
 \begin{align}
\xi^-  \in \omega_1^{\langle+ \rangle}(q_{(- m \times \per(q), 0]}).
 \end{align}

The statement $(-)$ results from (5.1) and (5.2).
 The case $(+)$ is symmetric.
\end{proof}

We adopt as a condition to be imposed on a hyposynchronizing subshift 
$X \subset \Sigma^\Bbb Z$ that satisfies $\bold H13$ - $\bold 19$ the following Condition $\bold H20$. In the formulation of this condition we use for an edge $e^- \in \mathcal E^-(X)$
($e^- \in \mathcal E^-(X)$), for $q,r \in P^0(X)$ and for
 $
 x \in  F^-(q,r)
 $
 ($
 x \in  F^+(r,q)
 $)
 the notation $\mu_{e^-}$ ($\mu_{e^+}$) for the number of occurences of the edge 
 $\mu_{e^-}$($\mu_{e^+}$) in the path $f^-(x)$($f^+(x)$).

\smallskip
 \noindent
$\bold H20$. 
 There exists $D\in \Bbb N$ such that the following statements $(-) $ and $(+)$ hold:

 $ (-)$ For $\xi^- \in\Xi^- (X)$ there exists a pair 
 $(e^- , \lambda ) \in \mathcal E^-\times \Bbb N$ such that one has for $q, r \in P^0(X)$, a sequence 
 $$
x^{(k)} \in F^-(q, r), \quad k \in \Bbb  Z_+,
 $$
 such that
 $$
 I^+( x^{(0)}) = 0, I^-( x^{(0)}) < -2D,
 $$
 $$
  I^+( x^{(k)}) = 0, I^-( x^{(k)}) =  I^-( x^{(0)}) - k,\quad (k \in \Bbb N)
 $$
 $$
 x^{(k)}_i = \xi^-,    \quad   I^-( x^{(k)}) + D \leq i \leq -D, \quad (k \in \Bbb Z_+)
 $$
 $$
  x^{(k)}_{(-\infty,  I^-( x^{(0)}) - k + D  )} =  x^{(0)}_{(-\infty,  I^-( x^{(0)})  + D  )},
  \quad (k \in \Bbb N)
 $$
and such that
  \begin{align}
 \lim_{k \to \infty}  x^{(k)} \in A^-(X),
  \end{align}
\begin{align}
& \ell (f^-(x^{(k)})=  \ell (f^-(x^{(0)})  +k\lambda, \quad k  \in \Bbb N,
 \\
& \mu_{e^-}(f^-(x^{(k)})=  \mu_{e^-} (f^-(x^{(0)})  +k\lambda, \quad k  \in \Bbb N.
 \end{align}
 $ (+)$ For $\xi^+ \in\Xi^+ (X)$ there exists an edge 
 $e^+  \in \mathcal E^+$ such that one has for
$q, r \in P^0(X)$, a sequence 
 $$
x^{(k)} \in F^+(r, q), \quad k \in \Bbb  Z_+,
 $$
 such that
 $$
 I^-( x^{(0)}) = 0, I^+( x^{(0)}) > 2D,
 $$
 $$
  I^-( x^{(k)}) = 0, I^+( x^{(k)}) =  I^+( x^{(0)}) + k,\quad (k \in \Bbb N)
 $$
 $$
 x^{(k)}_i = \xi^+,    \quad     D \leq i \leq I^+( x^{(k)})-D, \quad (k \in \Bbb Z_+)
 $$
 $$
  x^{(k)}_{(-\infty,  I^-( x^{(0)}) - k + D  )} =  x^{(0)}_{(-\infty,  I^-( x^{(0)})  + D  )},
  \quad (k \in \Bbb N)
 $$
and such that
  \begin{align}
 \lim_{k \to \infty}  x^{(k)} \in A^+(X),
  \end{align}
\begin{align}
& \ell (f^+(x^{(k)})=  \ell (f^+(x^{(0)})  +k, \quad k  \in \Bbb N,
 \\
& \mu_{e^+}(f^+(x^{(k)})=  \mu_{e^+} (f^+(x^{(0)})  +k, \quad k  \in \Bbb N.
 \end{align}

It follows from (5.6) and (5.7) that for $\xi^- \in \Xi^-(X)$ the pair $(e^-, \lambda)$ such that there exists a sequence $x^{(k)}, k\in \Bbb Z,$ as in  $\bold H20(-)$ is unique. Denoting this unique pair by 
$(e^-(\xi^-) , \lambda (\xi^-))$ we set 
$$
\Xi^-(e^-, \lambda) = \{\xi^-: (e^-(\xi^-) , \lambda (\xi^-))   = (e^-, \lambda) \} \quad \quad ((e^-, \lambda) \in \mathcal E^- \times \Bbb N),
$$

 Set
 $$
 \mu^-_X(e^-,\lambda) = \card (\Xi^-(e^-,\lambda)), \quad (e^-,\lambda) \in \mathcal E^-(X)  \times \Bbb N.
 $$
 We name the symbols in $\ \Xi^-(e^-,\lambda)$, enumerating them in the process,
 $$
 \Xi^-(e^-,\lambda) = \{ \xi^-(e^-,\lambda, \mu_\circ^-): 
  1\leq \mu_\circ^- \leq \mu_X^-(e^-,\lambda)   \},
 \quad (e^-,\lambda) \in \mathcal E^-(X)  \times \Bbb N.
 $$
 We set
 $$
 \Sigma^-(X) = \{ (e^-,\lambda, \mu_\circ^-):  
 (e^-,\lambda) \in \mathcal E^-(X)  \times \Bbb N,1 \leq  \mu_\circ^- \leq \mu_X^-(e^-,\lambda)    \}.
 $$

It follows from (5.9) and (5.10) that for $\xi^+ \in \Xi^-(X)$ the edge $e^+$ such that there exists a sequence $x^{(k)}, k\in \Bbb Z,$ as in  $\bold H20(+)$ is unique. Denoting this edge by 
$e^+(\xi^+)$ we set 
 $$
\Xi^+(e^+) = \{\xi^+: e^-(\xi^+)   = e^+ \}, \quad (e^+ \in \mathcal E^+).
$$ 
 Set
 $$
 \mu^+_X(e^+) = \card (\Xi^+(e^+)), \quad e^+ \in \mathcal E^+(X).
 $$
 We name the symbols in $\Xi^+(e^+)$, enumerating them in the process,
 $$
 \Xi^+(e^+) = \{ \xi^+(e^+, \mu_\circ^+):  1\leq \mu_\circ^- \leq \mu_X^-(e^+)   \},
 \quad e^+ \in \mathcal E^+(X).
 $$
 We set
 $$
 \Sigma^+(X) = \{ (e^+, \mu_\circ^+):  
 (e^+ \in \mathcal E^+(X),1 \leq  \mu_\circ^+ \leq \mu_X^-(e^+)    \}.
 $$

We have obtained bijections
$$
\Psi^-: \Sigma^-(X) \to \Xi^-(X), \quad \Psi^+: \Sigma^+(X) \to \Xi^+(X),
$$
 by
 $$
 \Psi^-( e^-, \lambda, \mu_\circ^-) =  \xi^-(e^-,\lambda, \mu_\circ^-),  \quad \quad 
( e^-, \lambda, \mu_\circ^-) \in  \Sigma^-(X),
$$ 
 $$
 \Psi^+( e^-, \mu_\circ^+) =  \xi^-(e^+,\lambda, \mu_\circ^+),  \quad \quad 
( e^+, \mu_\circ^+) \in  \Sigma^+(X).
$$ 
 
  \bigskip

We continue under the assumption that the hyposynchronizing subshift $X\subset \Sigma^\Bbb Z$ satisfies  $\bold H13$ - $\bold H20$ and that
\begin{align*}
P_\bullet^0(X) \neq \emptyset.
\end{align*}           
 In the next lemma and in the proof of the next theorem we choose a neutral periodic point  
$
p\in P_\bullet^0(X).
$
 This choice  has the effect that one has for all vertices $V^{hypo}$ of 
 $G^{hypo}_{\textcircled{\tiny V}}(X)$ that
 $$
 \tau_{p_{[0, \per(p))}} (V^{hypo}) = V^{hypo}.
 $$
This means that every vertex of $G^{hypo}_{\textcircled{\tiny V}}(X)$ has attached to it a loop with label sequence 
$p_{[0, \per(p))}$. One can insert into a path in $G^{hypo}_{\textcircled{\tiny V}}(X)$ after every edge any number of transversals of the loop with label sequence 
$p_{[0, \per(p))}$ that is attached to the target vertex of the edge. We will also use the reverse of this procedure

\begin{lemma}
$(-)$
For $(e^-,\lambda) \in \mathcal E^-(X)  \times \Bbb N$ and for $\xi^-\in\Xi^-(e^-,\lambda)$ one has that
$$
p_{(-\infty, 0)}\xi^-p_{[0, \infty)} \in (e^-)^\lambda.
$$
$(+)$
For $e^+ \in \mathcal E^+(X)$  and for $\xi^+ \in\Xi^+(e^+)$ one has that
$$
p_{(-\infty, 0)}\xi^+p_{[0, \infty)} \in e^+.
$$
\end{lemma}
\begin{proof}
The lemma follows from Lemma (5.2) in conjunction with $\bold H20$ (5.3), (5.4), (5.5), 
for 
$(-)$
and  $\bold H20$ (5.8),(5.9),(5.10),
for $(+)$.
\end {proof}

  \begin{theorem}
 Let $X \subset \Sigma^\Bbb Z$ be a hypersynchronizing subshift that satisfies 
 $\bold H13$ - $\bold H20$ such that
  $$
P_\bullet^0(X) \neq \emptyset.
$$
  Let the 
 Artin-Mazur
  zeta function of $X$ be equal to the Artin-Mazur of zeta function of   
 $
 Y_{\mu^-(X), \mu^+(X)}(G_\mathcal R(\{\textcircled{\tiny V} \}, \mathcal E^-(X) \cup \mathcal E^+(X)  ).
 $
 Then  the mapping $\Psi_X = (\Psi_X^- , \Psi_X^+  )$ implements a literal conjugacy of 
  $
 Y_{\mu^-(X), \mu^+(X)}(G_\mathcal R(\{\textcircled{\tiny V} \}, \mathcal E^-(X) \cup \mathcal E^+ (X) )
 $
 onto $X$.
 \end{theorem}
 \begin{proof}
 Given a word
 $$
 (e_i)_{1\leq i \leq I}\in \mathcal L( Y_{\mu^-(X), \mu^+(X)}(G_\mathcal R(\{\textcircled{\tiny V} \}, \mathcal E^-(X) \cup \mathcal E^+(X)  ))), \quad I \in \Bbb N, 
 $$
 one has by Lemma 5.2 that the concatenation of the points
 $$
 p_{(-\infty, 0)}\Psi(e_i)p_{[0, \infty)}
 $$
 is in $X$. Removing from this concatenation the transversals of the loops that have the label sequence 
 $
 p_{[0, \per (p)))}
 $
 and that are attached to the vertices of $G_{hypo}(X)$. shows that the word 
 $(\Psi(e_i))_{1\leq i \leq I}$ is admissible for $X$.
 It follows that the map $\Psi(X)$ is a continuous injection of  
 $Y_{\mu^-(X), \mu^+(X)}(G_\mathcal R(\{\textcircled{\tiny V} \}, \mathcal E^-(X) \cup \mathcal E^+(X)  )$ into $X$. The hypothesis on the Artin-Mazur zeta function of the subshifts $X$ and
 $Y_{\mu^-, \mu^+}(G_\mathcal R(\{\textcircled{\tiny V} \}, \mathcal E^- \cup \mathcal E^+  )$ implies that
 \begin{align}
 P(X) \subset \psi(Y_{\mu^-(X), \mu^+(X)}(G_\mathcal R(\{\textcircled{\tiny V} \}, \mathcal E^-(X) \cup \mathcal E^+(X)  ) ).
 \end{align}
It follows by $\bold H08$ that
 $$
 X= \psi(Y_{\mu^-(X), \mu^+(X)}(G_\mathcal R(\{\textcircled{\tiny V} \}, \mathcal E^-(X) \cup \mathcal E^+(X)  ) ) ). \qed
 $$ 
 \renewcommand{\qedsymbol}{}
  \end{proof}
  
  Renaming the symbols of its alphabet does not alter the structure of a subshift. It is therefore correct to say that  Theorem 5.1 amounts to a chacterization of the subshifts in the family $\mathcal Y$.
 To emphasize an aspect of the theorem we add a corollary.
\begin{corollary}
Let $X\subset \Sigma^\Bbb Z$ be a hyposynchronizing subshift, that satisfies
$\bold H13$ - $\bold 20$.
 Let the 
 Artin-Mazur
  zeta function of $X$ be equal to the Artin-Mazur of zeta function of   
 $
 Y_{\mu^-(X), \mu^+(X)}(G_\mathcal R(\{\textcircled{\tiny V} \}, \mathcal E^-(X) \cup \mathcal E^+(X)  ).
 $

Let the topological conjugacy class of $X$ contain a subshift 
$\bar X \subset \bar  \Sigma^\Bbb Z$ such that
$$
P_\bullet^0(\bar X) \neq \emptyset.
$$
 Then $\bar X$   is topologically conjugate to  
 $
 Y_{\mu^-(X), \mu^+(X)}(G_\mathcal R(\{\textcircled{\tiny V} \}, \mathcal E^-(X) \cup \mathcal E^+ (X) ).
 $
 \end{corollary}
We also add a corollary for the Dyck shift.    
 \begin{corollary}
 Let $X\subset \Sigma^\Bbb Z$ be a hyposynchronizing subshift, that satisfies
  $\bold H13$ - $\bold 20$.
 Let $N > 1$, and let
 $$
 \card (\mathcal E^-(X)) = \card (\mathcal E^+(X)) =N,
 $$.
 $$
 \card ( (\{ e^+ \in \mathcal E^+(X): (e^-,e^+) \in \mathcal R(X)\}) = 1, \qquad 
 (e^- \in \mathcal E^- (X)),
 $$
 $$
  \card ( (\{ e^- \in \mathcal E^-(X): (e^-,e^+) \in \mathcal R(X)\}) = 1, \qquad
  (e^+\in \mathcal E^+(X)).
 $$
 Let 
$$
\zeta(X) = \frac {2(1+\sqrt{1-4Nz^2})}{(1-2Nz+ \sqrt{1-4Nz^2})^2}.
$$
Let
$$
P_\bullet^0(X) \neq \emptyset.
$$
Then there is a literal conjugacy of $X$ onto the Dyck shift $D_N$.
  \end{corollary}
  
 \begin{proof}
 For the Artin-Mazur zeta function of the Dyck shift $D_N$ see \cite [Example 3]{Ke1}.
 \end{proof}

 \section{Examples}

 As starting point for the construction of examples we take the Dyck shift $D_2$. We denote by $\mathcal D_2$ the code that contains the words
 $(e_j)_{1 \leq j < I}$ that are defined by the conditions
 $$
 \prod_{1 \leq j\leq J} e_j   \not = \bold 1,  \quad  1\leq J < I,  
 $$
 $$
 \prod_{1 \leq \leq I} e_j   = \bold 1.
 $$
 \subsection{}
 One has a supply of hyposynchronizing subshifts that are obtained by excluding  correctly chosen words from $D_2$(see \cite{IK1}).
 
  \subsubsection{}
  Excluding from $D_2$ the word  $\alpha^-\alpha^+$ yields a hyposynchronizing subshift $X$ such that $P^0_\bullet(X) \neq \emptyset$.
   \subsubsection{}
   Excluding from $D_2$ the words  $\alpha^+\beta^+$and $\beta^+\alpha^+$ yields a hyposynchronizing subshift $X$ such that
$$
 \bigcap_{V^{hypo} \in \frak V^{hypo}(X)}   P^0(V^{hypo} )= \emptyset.
$$   
 
\subsection{}
We consider subshifts that appear as images of $D_2$ under conjugacies.

For a hyposynchronizing shift $X \in \Sigma^\Bbb Z$ that satisfies $\bold H12$ - 
$\bold H19$  we denote
 by $P^{0,0}(X)$ the set of $p \in \textcircled{\tiny V}$ such that $\vrt(p)$ has zero height. 
For hyposynchronizing subshifts $X$ that satisfy $\bold H13$ - $\bold H19$ one has by
$\bold H13$ that $P_\bullet(x) \neq \emptyset$ implies  
$\card (\{\vrt(p): p \in P^{0,0}(X)\}) = 1.$ This leads to. the observation that for hyposynchronizing subschifts the validity of $P_\bullet(x) \neq \emptyset$ is not invariant under topological conjugacy (Example 6.2.1. and Example 6.2.2.).

 \subsubsection{}
 Denote by $\widehat D$ the 2-block system of $D_2$ with alphabet
 $$
 \widehat \Sigma =
 \{(\alpha^-, \alpha^-), (\alpha^-, \alpha^+),  (\alpha^-, \beta^- ),(\alpha^+ , \alpha^-  ),(\alpha^+ , \alpha^+  ),(\alpha^+ , \beta^-  ), (\alpha^+ ,  \beta^+)\},
 $$
 and denote by $\varphi $ the conjugacy of $D_2$ onto $\widehat D$ that maps the point 
 $(x_i)_{i \in \Bbb Z} \in D_2$ to the point $(x_{i-1},x_i)_{i \in \Bbb Z} \in 
 \widehat \Sigma^\Bbb Z $. Denote by $\Phi$ the block map with memory one and anticipation zero that implements the conjugacy $\varphi$.
 
It holds that
$$
P^{0,0}(\widehat D) = \varphi(P^{0,0}(D_2)).
$$
We consider points  $ \widehat p \in P^{0,0}(\widehat D) $ distingishing the two cases $(\alpha) $ and $(\beta)$:
$$
(\alpha) \quad \quad \quad \widehat p_1 \in \{(\alpha ^+ , \alpha ^- ),  ( \alpha ^+, \beta^-  )\},
$$
$$
(\beta) \quad \quad \quad \widehat p_1 \in \{(\beta^+ , \alpha ^-  ),  (\beta^+ , \beta^-  )\}.
$$
In case $\alpha$ (case $\beta$) the  set  
$\omega^{\langle+ \rangle}(\widehat p_{(- \infty, 0]}) $  contains the points 
$\widehat{u} \in \widehat{\Sigma}^\Bbb N$ that are obtained by applying $\Phi$ to points 
$(u_i)_{0\leq i < \infty}$ with the initial symbol $u_0$ equal to $\alpha^+$($\beta^+$) that is followed by a right infinite concatenation of words in $\mathcal D_2$. Therefore
$$
\card (\{\vrt(p): p \in P^{0,0}(X)\}) = 2.
$$
\subsubsection{}
Let 
$$
Y \subset \{\alpha^- ,\alpha ^+ , \beta^- ,\beta ^+\}^\Bbb Z,
$$
be the image of $D_2$ under the conjugacy that is implemented by a block map with memory one and anticipation zero that interchanges the symbol $\alpha^-$ with the symbol $\alpha^+$ and the symbol $\beta^-$  with the symbol $\beta^+$ provided that these symbols are preceded by the symbol $\alpha^-$ or by the symbol $\alpha^+$:
\begin{align*}
&\Phi(\alpha^-,\alpha^-) =\Phi(\alpha^+ ,\alpha^-) =\Phi(\beta^+ , \alpha^+ ) = \alpha^-,
\\
&\Phi(\alpha^-,\alpha^+) =\Phi(\alpha^+,\alpha^+) =\Phi(\beta^- ,\alpha^-  ) = 
\Phi(\beta^+, \alpha^-  )=\alpha^+,
\\
&\Phi(\alpha^-, \beta^- ) =\Phi(\beta^-, \beta^- ) =\Phi(\beta^- , \beta^- ) = \Phi(\beta^- ,\beta^-  )=\beta^-,
\\
&\Phi(\alpha^+, \beta^+ ) =\Phi(\beta^-, \beta^-)  
=\Phi(\beta^-,  \beta^- ) = \beta^+.
\end{align*}

We consider a point
$$
y = \varphi(x), \quad  \quad\quad (x \in P^{0,0}(D_2))
$$
distinguishing the cases $\alpha$ and $\beta$:
$$
(\alpha) \quad  \quad\quad y_0 = \alpha^+,
$$
$$
(\beta) \quad  \quad\quad y_0 = \beta^+.
$$
In case $(\alpha)$(case $(\beta)$)
the set
$\omega^{\langle+ \rangle}(y_{(- \infty, 0]})$ contains the points in
$\{\alpha^- ,\alpha ^+ , \beta^- ,\beta ^+\}^\Bbb N$ that are obtained by applying $\Phi$ to points
$(u_i)_{0\leq i < \infty}$ that have an initial symbol in $\{\alpha^+,\beta^+\}$, followed by a right infinite concatenation of words in 
$\mathcal D_2$. In case ($\alpha$) (case ($\beta$)) the points in 
$\omega^{\langle+ \rangle}(y_{(- \infty, 0]})$ have an initial symbol in 
$\{\alpha^+, \beta^+  \}$($\{\alpha^-, \beta^- \}$). Therefore
$$
\card (\{\vrt(p): p \in P^{0,0}(X)\}) = 2.
$$

\subsubsection{}
Let 
$$
Y \subset \{\alpha^- ,\alpha ^+ , \beta^- ,\beta ^+\}^\Bbb Z,
$$
be the image of $D_2$ under the conjugacy that is implemented by a block map with memory one and anticipation zero that interchanges the symbol $\alpha^-$ with the symbol $\beta^-$ and that interchanges the symbol $\alpha^+$ with the symbol $\beta^-$ provided that these symbols are preceded by the symbol $\alpha^-$ or by the symbol 
$\beta^-$:
\begin{align*}
&\Phi(\alpha^-, \alpha^- ) =\Phi(\alpha^+ ,  \beta^-) =
\Phi(\beta^- , \alpha^- ) = 
\Phi(\beta^+ ,\beta^-  )=\beta^-,
\\
&\Phi(\alpha^+ , \alpha^- ) =\Phi(\alpha^- , \beta^- ) =\Phi(\beta^+ , \alpha^- ) = 
\Phi(\beta^-,\beta^-) =\alpha^-,
\\
&\Phi(\alpha^- ,\alpha^+  ) =\Phi(\alpha^+ , \beta^+ ) =\Phi(\beta^+ , \beta^+ ) = \beta^+,
\\
&\Phi(\alpha^+, \alpha^+ ) =\Phi(\beta^+ , \alpha^+ ) =
\Phi(\beta^- , \beta^+ ) = \alpha^+,
\end{align*}

It holds that
$$
P^{0,0}(\widehat Y) = \varphi(P^{0,0}(D_2)).
$$
For $y\in P^{0,0}(\widehat Y)$ one has that
$$
y_0 \in \{\alpha^+ , \beta^+\}.
$$
It follows for $x \in D_2$ and $y = \varphi(x)$ that one obtains the points in
$\omega^{\langle+ \rangle}(y_{(- \infty, 0]})$
by applying $\Phi$ to the points in $\omega^{\langle+ \rangle}(x_{(- \infty, 0]})$. 
Therefore
$$
\card (\{\vrt(p): p \in P^{0,0}(X)\}) = 1.
$$

\medskip

\par\noindent Wolfgang Krieger
\par\noindent Institute for Mathematics 
\par\noindent  University of Heidelberg
\par\noindent Im Neuenheimer Feld 205 
 \par\noindent 69120 Heidelberg
 \par\noindent Germany
\par\noindent krieger@math.uni-heidelberg.de

 \end{document}